\documentclass{article}
\usepackage[utf8]{inputenc}
\usepackage{cite}
\usepackage{hyperref}

\title{\vspace{-1.5cm}Niemeier lattices, smooth 4-manifolds and instantons}
\author{Christopher Scaduto\thanks{The author was supported by NSF grant DMS-1503100}\vspace{.1cm} \\ \vspace{.1cm} {\emph{Simons Center for Geometry and Physics}}\\ \texttt{cscaduto@scgp.stonybrook.edu}} 

\date{}

\usepackage[a4paper, total={6in, 8.5in}]{geometry}
\usepackage{graphicx}
\usepackage{amssymb}
\usepackage{tikz-cd}
\usepackage{titling}
\usepackage{mathrsfs} 
\usepackage{booktabs}
\usepackage[all]{xy}
\usepackage{amsthm}
\usepackage{diagbox}
\usepackage{tabularx}
\usepackage{amscd}
\usepackage{caption}
\usepackage{nicefrac}
\usepackage{amsmath}
\usepackage{mathtools}
\usepackage{tikz}
\usepackage{mathabx}
\usepackage{dsfont}
\usepackage{lipsum}
\usepackage{mwe}
\usepackage{slashed}
\usepackage{rotating}
\usepackage{subcaption}

\usepackage{xcolor}

\usepackage{xcolor}

\newcolumntype{Y}{>{\centering\arraybackslash}X}

\usepackage{sectsty}

\sectionfont{\fontsize{12}{15}\selectfont}

\usetikzlibrary{shapes,snakes}
\usetikzlibrary{matrix,arrows,decorations.pathmorphing}
\usetikzlibrary{positioning}

\newcommand{\R}{\mathbb{R}}
\newcommand{\C}{\mathbb{C}}
\newcommand{\Z}{\mathbb{Z}}

\newcommand{\Q}{\mathbb{Q}}
\newtheorem{theorem}{Theorem}[section]
\newtheorem{prop}[theorem]{Proposition}
\newtheorem{lemma}[theorem]{Lemma}

\newtheorem{corollary}[theorem]{Corollary}
\newtheorem{remark}[theorem]{Remark}

\begin{document}

\maketitle 
\vspace{-0.5cm}

\begin{abstract}
We show that the set of even positive definite lattices that arise from smooth, simply-connected 4-manifolds bounded by a fixed homology 3-sphere can depend on more than the ranks of the lattices. We provide two homology 3-spheres with distinct sets of such lattices, each containing a distinct nonempty subset of the rank 24 Niemeier lattices.
\end{abstract}

\vspace{-.2cm}

\section{Introduction}\label{sec:intro}

Let $X$ be a smooth, compact, oriented 4-manifold. The intersection form of $X$ is the free abelian group $L_X=H_2(X;\Z)/\text{Tor}$ equipped with the symmetric bilinear form $L_X\otimes L_X\to \Z$ defined by the intersection of 2-cycles. A well-known result of Donaldson \cite{d-connections} says that if $X$ has no boundary, and $L_X$ is positive definite, then it is equivalent over the integers to a diagonal lattice $\langle 1 \rangle^n$.\\

In general, a {\emph{lattice}} is a free abelian group $L$ of finite rank equipped with a symmetric bilinear form. We write $x \cdot y$ for the pairing of $x,y\in L$. A lattice $L$ is {\emph{unimodular}} if it has a basis $\{e_i\}$ for which $|\det(e_i\cdot e_j)|=1$. If $X$ as above has an integer homology 3-sphere boundary $Y$, then $L_X$ is a unimodular lattice. A lattice $L$ is {\emph{even}} if $x\cdot x$ is an even integer for every $x\in L$. A definite lattice is {\emph{minimal}} if it has no vectors of absolute norm $1$.\\

Fix an integer homology 3-sphere $Y$. Write $\mathscr{L}(Y)$ for the set of isomorphism classes of minimal definite unimodular lattices $L$ such that there exists a smooth, compact, oriented 4-manifold $X$ without 2-torsion in $H_\ast(X;\Z)$, $\partial X=Y$, and $L_X=L\oplus \langle 1\rangle^n$ for some $n\geqslant 0$. Write $\mathscr{E}_i(Y)\subset \mathscr{L}(Y)$ for the subset of even rank $i$ lattices. That $\mathscr{E}_i(Y)$ is empty for $i$ large enough was proven by Fr\o yshov using Seiberg-Witten theory \cite{froyshov-monopole}, and also follows from Heegaard Floer theory \cite{os}.\\

The rank of an even positive definite unimodular lattice is an integral multiple of $8$. There is only one such lattice of rank 8, and two of rank 16. In rank 24, there are 24 such lattices, the {\emph{Niemeier lattices}}, see e.g. \cite[Ch.6]{conwaysloane}, the set of which we denote by $\mathscr{N}_{24}$. The number beyond rank 24 grows quickly, and those of rank 32 already make up more than a billion. Write $S_{1}^3(K)$ for $+1$ Dehn surgery on a knot $K\subset S^3$, and $T_{n,m}$ for the positive $(n,m)$ torus knot.\\

\begin{theorem}\label{thm:2,11} $\mathscr{E}_{24}(S_{1}^3(T_{2,11}))$ and $\mathscr{E}_{24}(S_{1}^3(T_{4,5}))$ are distinct and nonempty subsets of $\mathscr{N}_{24}$.
\end{theorem}

\vspace{.3cm}

\noindent This result contrasts with the constraints on $\mathscr{E}_{24}(Y)$ imposed by Seiberg-Witten and Heegaard Floer theory. In those contexts, only the ranks of even definite lattices bounded by a homology 3-sphere $Y$ are constrained, and thus either $\mathscr{E}_{24}(Y)$ is empty or no further information is obtained. The knots $T_{2,11}$ and $T_{4,5}$ in Theorem \ref{thm:2,11} may be replaced by other knots; see \S \ref{sec:otherknots}.\\

The {\emph{nonempty}} part of Theorem \ref{thm:2,11} is well-known. Indeed, recalling that $S_{1}^3(T_{n,m})$ is diffeomorphic to the orientation-reversal of the Brieskorn sphere $\Sigma(n,m,nm-1)$, the canonical negative definite plumbings for $\Sigma(2,11,21)$ and $\Sigma(4,5,19)$ are even and of rank 24. The rest of Theorem \ref{thm:2,11} is new. The {\emph{distinct}} part of the theorem is provided by obstructions from instanton Floer theory, following variations on a method of Fr\o yshov \cite{froyshov-inequality} and recent work of the author \cite{scaduto-forms}. In the latter reference, we used relations in the instanton cohomology ring of a surface times a circle, taken with mod 2 and mod 4 coefficients; here, we utilize relations with mod 8 coefficients. Specifically, we will show that the plumbed lattice for $\Sigma(4,5,19)$ cannot occur for $\Sigma(2,11,21)$.\\

To put the above result into context, we define the {\emph{geometric 4-genus}} $g_4(L)$ of a minimal positive definite unimodular lattice $L$ to be the minimum $g\geqslant 0$ such that there exists a smooth, closed, oriented 4-manifold $W$ with $b_2^-(W)=1$ and no 2-torsion in its homology, and an embedded connected orientable surface $\Sigma\subset W$ of genus $g$ and self-intersection $-1$ such that the orthogonal complement of $[\Sigma]\in L_W$ is isomorphic to $L\oplus \langle + 1 \rangle^n$ for some $n\geqslant 0$.\footnote{A more natural definition omits the torsion hypothesis, only included here because our methods require it.} For a knot $K\subset S^3$ we have

\begin{equation}
    L \in \mathscr{L}(S_{1}^3(K)) \;\; \Longrightarrow \;\; g_4(K) \geqslant g_4(L). \label{eq:4gineq}
\end{equation}

\vspace{.2cm}

\noindent To see this, let $X$ be as in the definition of $\mathscr{L}(Y)$, so that $\partial X=Y$ and $L_X=L\oplus \langle 1 \rangle^n$. We may assume $L$ is nonstandard, or else $g_4(L)=0$; then $L$ is positive definite. Form $W$ by gluing to $X$ a standard 2-handle cobordism associated to $K$, and then filling in the resulting 3-sphere boundary with a 4-handle. A slicing surface for $K$ along with a disk in the 2-handle cobordism forms a surface $\Sigma$ of genus $g_4(K)$ and self-intersection $-1$. Then $W$ has $b_2^-(W)=1$, no 2-torsion, and the complement lattice of $[\Sigma]$ is $L_X$, so (\ref{eq:4gineq}) follows. This is the construction of \cite[Cor.1]{froyshov-inequality}. We compute:\\

\begin{theorem}\label{thm:g4}
    $g_4(D_{24})=5$ and $g_4(A_{24})=6$.
\end{theorem}

\vspace{.3cm}

\noindent Here, $D_{24}$ is the canonical plumbed lattice for $-\Sigma(2,11,21)$, and $A_{24}$ is that for $-\Sigma(4,5,19)$. Recall that $g_4(T_{2,11})=5$ and $g_4(T_{4,5})=6$. From Theorem \ref{thm:g4} and implication (\ref{eq:4gineq}), the lattice $A_{24}$ is not a member of $\mathscr{E}_{24}(S_{1}^3(T_{2,11}))$. Thus Theorem \ref{thm:2,11} follows from Theorem \ref{thm:g4}.\\

Seiberg-Witten and Heegaard Floer theory imply that if $L\in \mathscr{N}_{24}$, then $g_4(L)\geqslant 5$, see \S \ref{sec:oddsandends}. This is an equality for the Niemeier lattice $D_{24}$, as follows from (\ref{eq:4gineq}) and the observation that $D_{24}$ is the canonical plumbed lattice for $-\Sigma(2,11,21)=S_1^3(T_{2,11})$. Thus $g_4(D_{24})=5$ is not new. However, the other equality in Theorem \ref{thm:g4} is new. We will also show that the Niemeier lattice $A_{12}^2$ cannot occur for $-\Sigma(2,11,21)$ by showing that it has $g_4\geqslant 6$.\\

\vspace{.25cm}

\textbf{Outline.} In Section \ref{sec:obstructions}, the obstructions derived from instanton Floer theory are introduced. These are applied in Section \ref{sec:niemeier} to several Niemeier lattices, where Theorem \ref{thm:g4} is proved. In Section \ref{sec:oddsandends}, we discuss some other applications and remaining questions.\\

\vspace{.2cm}

\textbf{Acknowledgments.} The author would like to thank Marco Golla for many inspiring conversations, Stefan Behrens for some comments, and Simon Donaldson for his encouragement and support. The author was supported by NSF grant DMS-1503100.\\

\vspace{.3cm}

\newpage

\section{The obstructions}\label{sec:obstructions}

In this section we review the obstructions used to obtain lower bounds on the geometric 4-genus of a given lattice. Most of this material runs parallel to \cite[\S 2]{scaduto-forms}. The difference is that we also consider inequalities derived from instanton Floer theory with coefficients modulo 8.  \\

Let $L$ be a positive definite unimodular lattice. We write $x\cdot y$ for the pairing of elements $x,y\in L$ and $x^2 = x\cdot x$. The {\emph{norm}} of $x\in L$ is defined to be $x^2$. For any subset $S\subset L$ denote by $\text{Min}(S)$ the set of elements of minimal norm among all elements within $S$. We say $w\in L$ is {\emph{extremal}} if it is of minimal norm within its index 2 coset, i.e. $w\in \text{Min}(w+2L)$. As in \cite{scaduto-forms}\footnote{Actually, our $f_2(L)$ here is a priori only less than or equal to that of \cite{scaduto-forms}; it is equal to $m(L)$ of that reference.}, define

\[
    f_2(L) \; := \; \max\left\{ w^2-1 : \; w\in L,\; w \text{ extremal},\; |\text{Min}(w+2L)/\pm| \not\equiv 0\; (\text{mod 2}) \right\}.
\]

\vspace{.2cm}

\noindent Following \cite{froyshov-equivariant}, for $w,a\in L$ and $m\geqslant 0$ such that $w^2 \equiv m$ (mod 2), define

\begin{equation}
    \eta(L,w,a,m) \; := \;  \frac{1}{2}\sum_{z\in\text{Min}(w+2L)}(-1)^{((z+w)/2)^2}(a\cdot w)^m.\label{eq:eta}
\end{equation}

\vspace{.2cm}

\noindent When $m=0$ we interpret $(a\cdot w)^m=1$ and simply write $\eta(L,w)$. When $w$ is of even norm and $w\not\equiv 0$ (mod 2), $\eta(L,w)$ is a signed count of the elements in $\text{Min}(w+2L)/\pm$. \noindent For integers $n>2$ set

\[
    f_n(L) \; := \; \max\left\{ \frac{w^2-m}{2}: \; \begin{array}{c}\exists\; w,a\in L, m\geqslant 0,\; w^2 \equiv m \,(\text{mod }2),\\  w \text{ extremal},\;\; \eta(L,w,a,m)\not\equiv 0 \,(\text{mod }n) \end{array}\right\}. 
\]

\vspace{.2cm} 

\noindent For a negative definite unimodular lattice $L$ we define all of the above quantities using $-L$. Let $\Sigma$ be a Riemann surface of genus $g$. Let $R$ be a ring. The instanton cohomology $I^\ast(\Sigma\times S^1;R)$ for the pull-back of a $U(2)$-bundle with odd first Chern class over $\Sigma$ has a ring structure, studied by Mu\~{n}oz \cite{munoz} when $R=\C$. There are distinguished elements $\alpha,\beta,\gamma\in I^\ast(\Sigma\times S^1;R)$ which in Donaldson-Floer theory correspond to $\mu$-classes for the relative invariants of $\Sigma\times D^2$. We use the conventions of \cite{scaduto-forms}. Define $\smash{N_\alpha^n(g)}$ and $\smash{N_\beta^n(g)}$ to be the nilpotency degrees of $\alpha$ and $\beta$, respectively, in the ring $I^\ast(\Sigma\times S^1; \Z)\otimes\Z/n$ modulo $\gamma$.\\

\begin{theorem}\label{thm:inequality}
    Let $X$ be a smooth, closed, oriented 4-manifold with all torsion in its homology coprime to $n\geqslant 2$. Suppose $b_2^+(X)=1$, and let $\Sigma\subset X$ be an embedded surface with $\Sigma\cdot \Sigma= 1$ and genus $g$. Define $L\subset L_X$ to be the negative definite unimodular lattice orthogonal to $[\Sigma]$.  Then
        \begin{align}
        N_\alpha^2(g) \;  & \geqslant \; f_2(L) \qquad (n =2) \label{eq:mainineq2} \\
        N_\beta^n(g) \; & \geqslant \; f_n(L)\qquad (n >2)\label{eq:mainineq}
    \end{align}
\end{theorem}

\vspace{.35cm}

\noindent This follows from Thms. 2.1 and 2.2 of \cite{scaduto-forms} when $n=2$ and $n=4$, respectively, but the proof therein of the latter, which is adapted from \cite{froyshov-inequality}, works for all integers $n > 2$.\\

In \cite[Prop.2.3]{scaduto-forms} we showed that $\smash{N_\alpha^2(g)=g}$ and $\smash{N_\beta^4(g)=\lceil g/2\rceil}$ for $g\leqslant 128$, and conjectured that these identities hold in general. To this we add the following, which, as mentioned in the introduction, will be used to prove the {\emph{distinct}} part of Theorem \ref{thm:2,11}:\\

\begin{prop}\label{prop:8}
    For $g\leqslant 128$ with $g\not\equiv 2$ \emph{(mod 4)}, we have $N_\beta^8(g)=\lceil g/2\rceil$.
\end{prop}

\vspace{.2cm}

\begin{proof}
    The proof is essentially the same as for the cases $n\in\{2,4\}$ from \cite[\S 6]{scaduto-forms}; for completeness we sketch the argument. Let $R=\Z[\alpha,\beta]$. We define  $\zeta_r=\zeta_r(\alpha,\beta)\in R$ recursively by $\zeta_0=1$, and

\[
    \zeta_{r+1} \; = \; \alpha \zeta_r + r^2(\beta+(-1)^r8)\zeta_{r-1}.
\]

\vspace{.2cm}

\noindent Define the ideals $J_r=(\zeta_r,\zeta_{r+1},\zeta_{r+2})\subset R$. Then $R/J_g\otimes \Q$ is isomorphic to the $\text{Sp}(2g,\Z)$-invariant $\mathbb{Z}/4$-graded instanton cohomology ring of $\Sigma\times S^1$ modulo $\gamma$, with rational coefficients. This was proven by Mu\~{n}oz \cite{munoz}. Indeed, the $\zeta_r$ we have defined are obtained from those in \cite[\S 6]{scaduto-forms} by setting $\gamma=0$. Here $\Sigma$ is a Riemann surface of genus $g$, and the $U(2)$-bundle used in defining the instanton cohomology is the pullback of a bundle over $\Sigma$ with odd first Chern class.\\

The elements $\alpha$ and $\beta$ correspond to integral generators in the instanton cohomology ring. Furthermore, according to Lemma 6.3 of \cite{scaduto-forms}, there is an integral generator $\varepsilon$ such that $\beta-\alpha^2=8\varepsilon$, so that $\beta\equiv \alpha^2$ (mod 8). We define the following polynomials:

\[
    \theta_r(\alpha,\varepsilon) \; := \; \frac{(2r-3)!!}{r!}\zeta_r(\alpha,\beta)|_{\beta=\alpha^2+8\varepsilon}
\]

\vspace{.2cm}

\noindent Here $k!!=k(k-2)(k-4)\cdots (3)(1)$ for odd positive $k$, and $k!!=1$ if $k<0$. These polynomials are obtained from the ones in \cite[Conj. 6.4]{scaduto-forms} by setting $\gamma=0$. There we conjectured, and verified for many $r$, that each $\theta_r(\alpha,\varepsilon)$ has integer coefficients, and that $\theta_r(\alpha,\varepsilon)$ reduces mod 4 to $\pm\alpha^r$. This implies $N_\alpha^2(g)\leqslant g$ and, as $\beta\equiv \alpha^2$ (mod 4), that $N_\beta^4(g)\leqslant \lceil g/2 \rceil$. The additional observation we make here is that if $r\equiv 0,1$ (mod 4), $\theta_r(\alpha,\varepsilon)$ also reduces mod 8 to $\pm\alpha^r$. For example, $\theta_0= 1$, $\theta_1= \alpha$, 
\[
\theta_4(\alpha,\varepsilon) \; =\;  15  \alpha^{4} + 160  \alpha^{2}  \varepsilon - 120  \alpha^{2} + 360  \varepsilon^{2} - 720  \varepsilon + 360\; \equiv \;-\alpha^4 \text{ (mod 8)},
\]
\[
\theta_5(\alpha,\varepsilon) \; =\;  105  \alpha^{5} + 1456  \alpha^{3}  \varepsilon + 840  \alpha^{3} + 4984  \alpha  \varepsilon^{2} + 6160  \alpha  \varepsilon + 3192  \alpha\; \equiv \;\alpha^5 \text{ (mod 8)}.
\]

\vspace{.2cm}

\noindent These cases are the only ones needed for the applications of the current article, but we have verified the claim for $r\leqslant 128$ using a computer. Therefore, in this range, for $g\equiv 0,1$ (mod 4) we have $\smash{N^8_\beta(g)\leqslant \lceil g/2\rceil}$. From the definitions, $N^8_\beta(g)\geqslant N^4_\beta(g)=\lceil g/2\rceil $, so equality holds for these cases. As $N_\alpha^2(g)$ and $N_\beta^n(g)$ are non-decreasing functions in $g$, see e.g. \cite[Cor.19]{munoz}, for $g\equiv 3$ (mod 4) we have $N^4_\beta(g)=\lceil g/2\rceil \leqslant\ N_\beta^8(g)\leqslant N_\beta^8(g+1)=\lceil (g+1)/2 \rceil=\lceil g/2\rceil $, implying the result for this case.
\end{proof}

\vspace{.35cm}

\noindent From Theorem \ref{thm:inequality}, Proposition \ref{prop:8} and \cite[Cor.3.1]{scaduto-forms}, we obtain the following.\\

\begin{corollary}\label{cor:lowerbounds} Let $L$ be a definite unimodular lattice and suppose $g_4(L)\leqslant 128$. Then
    \begin{itemize}
        \item[(i)] $g_4(L) \geqslant f_2(L)$.
        \item[(ii)] $g_4(L) \geqslant 2f_4(L)-1$.
        \item[(iii)] If $f_8(L)\equiv 0$ \emph{(mod 2)}, then $g_4(L)\geqslant 2f_8(L)-2$.
        \item[(iv)] If $f_8(L)\equiv 1$ \emph{(mod 2)}, then $g_4(L)\geqslant 2f_8(L)-1$.
    \end{itemize}
\end{corollary}

\vspace{.3cm}

\noindent 
\begin{remark}{\emph{
Inequalities derived from cases in which $n=2^k$ is a higher power of 2 may also be obtained, but we will not make use of these in this article. We mention, however, that we have verified for low values of $r$ the following: the polynomial $\smash{(2r-3)!!\alpha^{\sigma}\zeta_r/r! |_{\alpha^2 = \beta-8\varepsilon}}$ lies in $\Z[\beta,\varepsilon]$, where $\sigma=0$ (resp. 1) if $r$ is even (resp. odd), and it reduces to $-\beta^{(r+\sigma)/2}$ modulo $2^k$ for $r\equiv 0$ (mod $2^{k-1}$) and $k\in\{2,3,4,5,6\}$. This implies that $\smash{N_\beta^n(g)=\lceil g/2\rceil}$ for $g\equiv 0$ (mod $2^{k-1}$) and $k\in\{2,3,4,5,6\}$. Note for these values of $n$ that $\smash{N_\beta^n(g)\leqslant 2\lceil g/2\rceil}$ always holds, because $\smash{(\beta^2-64)^{\lceil g/2\rceil}}$ vanishes modulo $\gamma$ in the ring with integer coefficients, see \cite[\S 5]{munoz}.
}}
\end{remark}

\subsubsection{Another expression for $\eta(L,w)$}

\noindent When computing $f_n(L)$, it will be useful to rewrite $\eta(L,w)$. Fix a positive definite unimodular lattice $L$. Let $w\in L$ be extremal and of even norm. For $i\in \Z$ define

\[
    S_{i}^w \; := \; \left\{ u\in L: \; u^2=i=-u\cdot w \right\}.
\]

\vspace{.2cm}

\noindent Note $S_i^w$ is empty for $i<0$ and $S_0^w=\{0\}$. Let $z\in \text{Min}(w+2L)$ and write $z=w+2u$. Then $w^2=z^2=(w+2u)^2=w^2+4w\cdot u + 4u^2$ implies $u\cdot w = -u^2$, and thus $\smash{u\in S_i^w}$ where $i=u^2$. Further, $\smash{(-1)^{((z+w)/2)^2}=(-1)^{(w+u)^2}=(-1)^{i}}$. Thus the sum (\ref{eq:eta}) with $m=0$ may be rewritten as $\eta(L,w) = \frac{1}{2}\sum_{i \geqslant 0}(-1)^{i}|S_i^w|$. Next note that $u\mapsto -w-u$ is a bijection from $S_i^w$ to $S_{m-i}^w$ where $m=w^2$. Thus $|S_i^w|=|S_{m-i}^w|$ for all $i\in \Z$ and in particular $|S_i^w|=0$ for $i>m$. Finally, we have

\begin{equation}
    \eta(L,w) \; = \; 1 - |S_1^w| + |S_2^w| - |S_3^w| + \cdots + (-1)^{m/2-1} |S_{m/2-1}^w| + (-1)^{m/2} \frac{1}{2}|S_{m/2}^w|.\label{eq:etanew}
\end{equation}

\vspace{.2cm}

\noindent When $L$ is even and $w$ is of norm 6, we have $\eta(L,w)=1+|S_2^w|$; when $w$ is instead of norm 8, $\eta(L,w)=1+|S_2^w|+\frac{1}{2}|S_4^w|$; and when $w$ is of norm 10, $\eta(L,w)=1+|S_2^w|+|S_4^w|$.\\

\subsubsection{Extremality criterion} It is also useful to have some method to verify that $w\in L$ is extremal. In general, we have

\begin{equation}
\textstyle
w \text{ is extremal} \quad \Longleftrightarrow \quad \text{  for all }u\in L\text{ with }u^2< \frac{w^2}{2} \text{ we have } |u\cdot w|\leqslant u^2 \label{eq:ex}
\end{equation}

\vspace{.2cm}

\noindent To see this, first suppose $w$ is not extremal, i.e. there is $u\in L$ such that $z=w+2u$ satisfies $z^2<w^2$. Upon possibly replacing $z$ by its negative, we may assume $z\cdot w \geqslant 0$. This implies $4u^2= (z-w)^2 \leqslant  z^2 + w^2 < 2w^2$, and so $u^2<w^2/2$. Furthermore, $w^2 > z^2 = (w+2u)^2 = w^2 + 4w\cdot u +4u^2$ implies $|w\cdot u| > u^2$. The converse follows similarly. Criterion (\ref{eq:ex}) implies, for example, that when $L$ is even, to show $w\in L$ of norm 6 or 8 is extremal, we only need show $|u\cdot w|\leqslant 2$ for all $u\in L$ of norm 2.\\

\vspace{.5cm}

\section{Niemeier lattices}\label{sec:niemeier}

In this section we prove Theorem \ref{thm:g4}. For several Neiemeier lattices $L$, we bound $f_n(L)$ from below for some $n\in\{2,4,8\}$, which via Corollary \ref{cor:lowerbounds} provides a lower bound for $g_4(L)$. We then provide upper bounds by producing examples of 4-manifolds with embedded surfaces.\\

We now fix some notation and provide some preliminaries. We refer the reader to the monograph \cite{conwaysloane} for general background. The classical root lattices of type ADE are defined as follows:
\begin{align*}
    \textstyle
    \mathsf{A}_n \; & := \; \left\{  (x_i)\in \Z^{n+1}: \; \sum x_i =0\right\} \qquad &\det\mathsf{A}_n &=n+1 \vspace{.1cm}\\
    \textstyle
    \mathsf{D_n} \; &:= \; \left\{  (x_i)\in \Z^n: \; \sum x_i\equiv 0 \text{ (mod 2)}\right\} \qquad &\det \mathsf{D}_n &=4 \vspace{.1cm}\\
    \textstyle
    \mathsf{E}_6 \; &:= \;  \textstyle \left\{ (x_i)\in {\textstyle{\frac{1}{2}}}\Z^{8}: \; x_i-x_j\in \Z,\;\sum_{i=1}^6x_i =x_7+x_8=0 \right\} \qquad &\det\mathsf{E}_6 &= 3 \vspace{.1cm}\\
    \textstyle
    \mathsf{E}_7 \; &:= \;  \left\{ (x_i)\in {\textstyle{\frac{1}{2}}}\Z^{8}: \; x_i-x_j\in \Z,\;\sum x_i =  0 \right\} \qquad &\det\mathsf{E}_7 &= 2 \vspace{.1cm}\\
    \textstyle
    \mathsf{E}_8 \; &:= \;  \left\{ (x_i)\in {\textstyle{\frac{1}{2}}}\Z^{8}: \;x_i-x_j\in \Z,\; \sum x_i\equiv 0 \text{ (mod 2)} \right\} \qquad &\det\mathsf{E}_8 &= 1 \vspace{.1cm}
\end{align*}

\noindent In general, a {\emph{root lattice}} is a positive definite lattice generated by vectors of norm 2, called \emph{roots}. It is well-known that every root lattice is some direct sum of the above listed lattices. For a given positive definite unimodular lattice $L$ with no vectors of norm 1, the {\emph{root lattice of $L$}} is the sublattice of $L$ generated by roots. \\

The 24 Niemeier lattices are listed in \cite[Ch.16]{conwaysloane}, labelled by their root lattices. For example, the Niemeier lattice $D_{24}$ is the unique even positive definite unimodular lattice of rank 24 with root lattice $\mathsf{D}_{24}$. Note that we have used different fonts for root lattices and unimodular lattices. Of the 24 Niemeier lattices, 23 have full rank root lattices, and one has no roots: the Leech lattice $\Lambda_{24}$.\\

For a lattice $L$ define the dual lattice $L^\ast=\{x\in L\otimes \Q: x\cdot y\in \Z\; \forall\, y\in L\}$. The determinant of $L$ is equal to  the order of the discriminant group $L^\ast/L$. Each Niemeier lattice $L$ apart from the Leech lattice is generated by its root lattice and a finite number of {\emph{glue vectors}}, cf. \cite[Ch.4 \S 3]{conwaysloane}. Suppose the root lattice of $L$ is $R_1\oplus \cdots \oplus  R_m$, with each $R_i$ a root lattice of type ADE. The glue vectors for $L$ are then elements of $R_1^\ast\oplus \cdots \oplus R^\ast_m$ that descend to generate a subgroup of the discriminant group $R_1^\ast/R_1\oplus \cdots \oplus R^\ast_m/R_m$ of square-root index. Glue vectors for Niemeier lattices can be read off from Table 16.1 in \cite[Ch.16]{conwaysloane}.\\

The {\emph{theta series}} of a positive definite lattice $L$ is the power series $\theta_L(q)=\sum a_i q^{i}$ where $a_i=a_i^L$ is the number of vectors in $L$ of norm $i$. By the classical theory of modular forms, the theta series of a Niemeier lattice $L$ is a linear combination of $\theta_{E_8^3}(q)$ and $\theta_{\Lambda_{24}}(q)$. We have
\begin{gather*}
    \theta_{E_8^3}(q) \; = \; 1 + 720q^2 + 179280q^4  + 16954560q^6 +  396974160q^8 + \ldots\\
    \theta_{\Lambda_{24}}(q) \; = \; 1 +  196560q^4 +  16773120q^6 +  398034000q^8 + \ldots
\end{gather*}

\vspace{.2cm}

\noindent From $\theta_L(q)$ having constant coefficient 1, we deduce $\theta_L(q)=\frac{a_2}{720}\theta_{E_8^3}(q) + (1 - \frac{a_2}{720})\theta_{\Lambda_{24}}(q)$. Thus all coefficients $a_i$ may be written in terms of $a_2$. For example, 
\begin{gather}
    a_4 \; = \; 196560 - 24a_2  \label{eq:a4}\\
    a_6 \; = \; 16773120 + 252a_2 \nonumber\\
    a_8 \; = \; 398034000 - 1472 a_2\nonumber
\end{gather}

\vspace{.35cm}

\textbf{\color{red}Notation:} We frequently write superscripts for repeated entries, i.e. $(\nicefrac{3}{2},\nicefrac{1}{2}^2)=(\nicefrac{3}{2},\nicefrac{1}{2},\nicefrac{1}{2})$. When we say {\emph{a vector of type}} $(\nicefrac{3}{2},\pm\nicefrac{1}{2}^2)$, we mean a vector obtained from $(\nicefrac{3}{2},\nicefrac{1}{2}^2)$ by permuting coordinates and changing signs where indicated. Thus there are 12 vectors of type $(\nicefrac{3}{2},\pm\nicefrac{1}{2}^2)$.\\

%The notation ${}_e\pm$ (resp. ${}_o\pm$) indicates that only even (resp. odd) numbers of negative signs are considered. For example, there are 2 vectors of type $({}_o\hspace{-.08cm}\pm \nicefrac{1}{2}^2)$, given by $(-\nicefrac{1}{2},\nicefrac{1}{2})$ and $(\nicefrac{1}{2},-\nicefrac{1}{2})$.

\vspace{.2cm}

\subsection{The lattice $D_{24}$}

Here we compute $g_4(D_{24})=5$, the first part of Theorem \ref{thm:g4}. As mentioned in the introduction, this is not new, but we include an argument for completeness, and in preparation for later arguments.\\

The Niemeier lattice $L=D_{24}$ is generated by the root lattice $\mathsf{D}_{24}$ and the glue vector $g=(\nicefrac{1}{2}^{24})$. Note $g^2=6$. It is easy to verify that $w=g$ is extremal. For example, the roots in $\mathsf{D}_{24}$ are vectors of type $(\pm 1^2, 0^{22})$; any such root $u$ has $|u\cdot w|\in \{0,1\}$, and so (\ref{eq:ex}) is satisfied. Furthermore, $\text{Min}(w+2L)=\{w,-w\}$. Indeed,  as just noted, $S^w_{2}=\emptyset$, and thus $\eta(L,w)=1+|S_2^w|=1$. It follows that $f_2(L)\geqslant w^2-1 = 5$, and by Corollary \ref{cor:lowerbounds} (i) that $g_4(L)\geqslant 5$. This was in fact already computed in \cite[Prop.4.1]{scaduto-forms}, and $D_{24}$ belongs to the family of lattices $\Gamma_{4k}$ from \cite[Prop.1]{froyshov-inequality}.\\

On the other hand, $D_{24}$ is isomorphic to the canonical plumbed lattice for $-\Sigma(2,11,21)$, given as follows, in which each unmarked node has weight 2:
\begin{center}
\begin{tikzpicture}[scale=.5]

	\draw (0,0) -- (22,0);
	\draw (2,0) -- (2,1);
	
	\draw[fill=black] (0,0) circle(.1);
	\draw[fill=black] (1,0) circle(.1);
	\draw[fill=black] (2,0) circle(.1);
	\draw[fill=black] (3,0) circle(.1);
	\draw[fill=black] (4,0) circle(.1);
	\draw[fill=black] (5,0) circle(.1);
	\draw[fill=black] (6,0) circle(.1);
	\draw[fill=black] (7,0) circle(.1);
	\draw[fill=black] (8,0) circle(.1);
	\draw[fill=black] (9,0) circle(.1);
	\draw[fill=black] (10,0) circle(.1);
	\draw[fill=black] (11,0) circle(.1);
	\draw[fill=black] (12,0) circle(.1);
	\draw[fill=black] (13,0) circle(.1);
	\draw[fill=black] (14,0) circle(.1);
	\draw[fill=black] (15,0) circle(.1);
	\draw[fill=black] (16,0) circle(.1);
	\draw[fill=black] (17,0) circle(.1);
	\draw[fill=black] (18,0) circle(.1);
	\draw[fill=black] (19,0) circle(.1);
	\draw[fill=black] (20,0) circle(.1);
	\draw[fill=black] (21,0) circle(.1);
	\draw[fill=black] (22,0) circle(.1);
	\draw[fill=black] (2,1) circle(.1);
	
	\node at (-.5,.25) {$6$};
	 
\end{tikzpicture}
\end{center}
Indeed, the node of weight 6 corresponds to $g$, and the nodes of weight 2 correspond to $(1,1,0^{22})$, $(1,-1,0^{22})$, $(0,1,-1,0^{21})$, \ldots, $(0^{21},1,-1,0)$. As $-\Sigma(2,11,21)$ is $+1$ surgery on $T_{2,11}$ of slice genus 5, by inequality (\ref{eq:4gineq}) we have $g_4(D_{24})\leqslant 5$, implying equality.\\

\subsubsection{From a collection of lines}

\begin{figure}
\centering
\begin{minipage}{.5\textwidth}
  \centering
  \includegraphics[scale=.4]{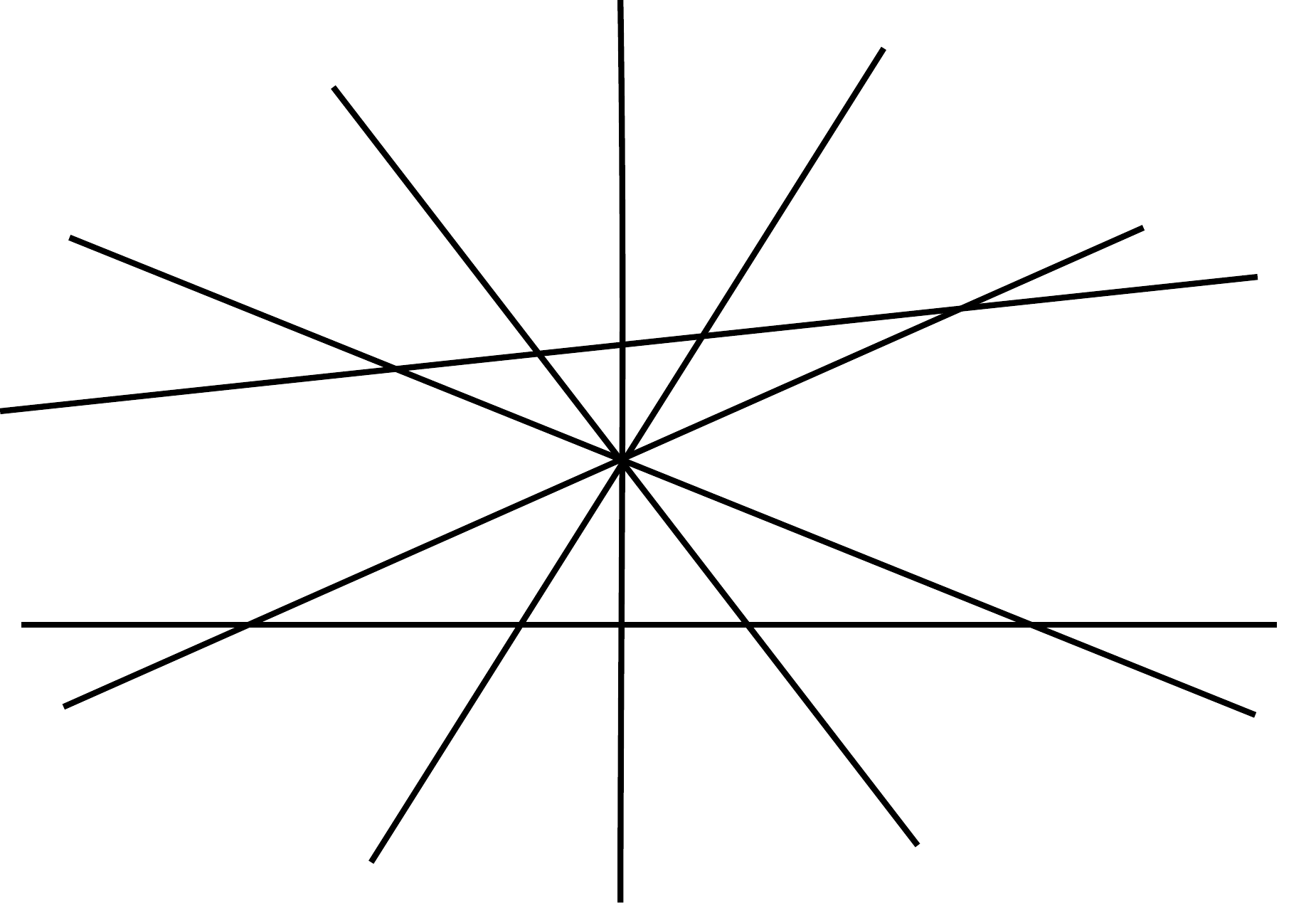}\vspace{.15cm}
  \captionof{figure}{$D_{24}$}
  \label{fig:d24}
\end{minipage}%
\begin{minipage}{.5\textwidth}
  \centering
  \includegraphics[scale=.4]{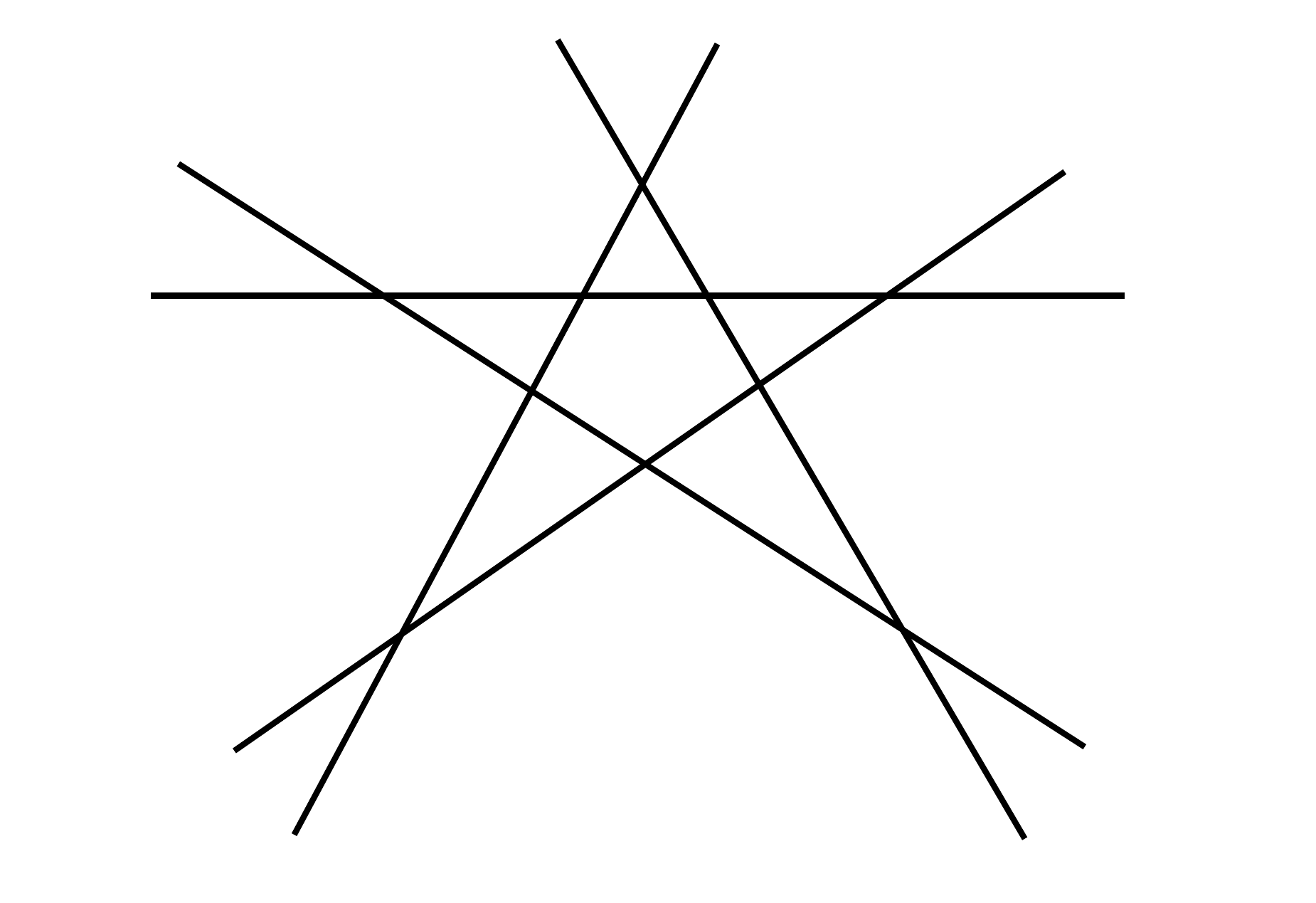}\vspace{.15cm}
  \captionof{figure}{$A_{24}$}
  \label{fig:a24}
\end{minipage}
\end{figure}

We provide another argument for $g_4(D_{24})\leqslant 5$, of a sort we return to for later use. Take a collection of 7 lines in the complex projective plane $\mathbb{C}\mathbb{P}^2$ that has 5 lines intersecting in one point, and no other points of multiplicity greater than 2, as in Figure \ref{fig:d24}. Blow up at the point of multiplicity 5, and 23 generic points. We slightly perturb the resulting curve to remove all double points, cf. \cite[p.39]{gs}. Then we obtain a smooth complex curve $\Sigma$ in $\smash{\mathbb{C}\mathbb{P}^2\# 24\overline{\mathbb{C}\mathbb{P}}^2}$ with

\begin{equation}
    [\Sigma] \; = \; 7h - 5e_1 -e_2 -\cdots -e_{24} \; = \; (7 \, |\, 5,1^{23}).\label{eq:firstsurface}
\end{equation}

\vspace{.2cm}

\noindent Note that $\Sigma$ is connected. Indeed, we may get from any one line to another in the original arrangement by passing through only double points; this property holds for similar constructions in the sequel, as we leave the reader to check. In (\ref{eq:firstsurface}), $h$ is induced from a generator for $H_2(\mathbb{C}\mathbb{P}^2;\Z)$, and each $e_i$ is an exceptional sphere associated to the $i^\text{th}$ blow-up. Note $[\Sigma]^2=1$. The notation $(7\, |\, 5,1^{23})$ follows \cite{conwaysloane}; in general, we write $(a \, |\, b_1,\ldots ,b_n)$ for an element $a h - \sum{b_i}e_i$ in the Lorentzian lattice $H_2(\smash{\mathbb{C}\mathbb{P}^2\# n\overline{\mathbb{C}\mathbb{P}}^2};\Z)$, and use superscripts for repeated entries. The canonical class $K$ of $\smash{\mathbb{C}\mathbb{P}^2\# 24\overline{\mathbb{C}\mathbb{P}}^2}$ is Poincar\'{e} dual to $\sum e_i-3h$, and the adjunction formula, in general given by

\begin{equation}
    0 \; = \; \langle K, [\Sigma]\rangle + [\Sigma]^2 + 2 - 2g(\Sigma),\label{eq:adjunction}
\end{equation}

\vspace{.2cm}

\noindent in this case yields that the genus of $\Sigma$ is 5. The lattice $L$ orthogonal to $[\Sigma]$ has roots $e_2-e_3,\ldots,e_{23}-e_{24}$, $h-e_1-e_2-e_3$ and $6h-4e_1-e_2-\cdots-e_{23}$, generating a copy of $-\mathsf{D}_{24}$. Furthermore, any orthogonal complement defined by a vector $v=(a \, | \, b_1,\ldots, b_n)$ with $a$ and all $b_i$ odd is an even lattice, as $x\cdot v=0$ implies $x^2 = \sum x_i^2 \equiv \sum x_i  \equiv ax_0+\sum b_i x_i = 0$ (mod 2). Thus $L$ is isomorphic to $-D_{24}$, and we conclude again that $g_4(D_{24})\leqslant 5$.\\

\subsection{The lattice $A_{24}$}

In this section we compute $g_4(A_{24})=6$, the second part of Theorem \ref{thm:g4}.\\

The lattice $L=A_{24}$ is generated by its roots and the glue vector $g=(\nicefrac{1}{5}^{20}, -\nicefrac{4}{5}^{5})$, which has $g^2=4$. Consider $w=(1^4,-1^4,0^{17})$ of norm 8. All roots $u\in L$ are of type $(1,-1,0^{23})$, and so satisfy $|u\cdot w|\leqslant 2$. Thus (\ref{eq:ex}) is satisfied and $w$ is extremal. We now compute $\eta(L,w)$. First, there are 16 roots such that $u\cdot w=-2$, namely $u=(-1,0^3,1,0^{20})$ and its permutations preserving coordinates 1 to 4, and 5 to 8. Thus $|S_2^w|=16$. To compute $|S_4^w|$, we list vectors of norm 4 in $A_{24}$ by type:

\begin{itemize}
    \item[(i)] $\phantom{\pm}(1^2,-1^2,0^{21})$
    \item[(ii)] $\pm (\nicefrac{1}{5}^{20}, -\nicefrac{4}{5}^{5})$
\end{itemize}

\noindent There are ${4\choose 2}{25\choose 4}=75900$ vectors of type (i), and $2{25 \choose 5}=106260$ of type (ii), summing to $182160$. This agrees with the theta series coefficient $a_4$ computed from (\ref{eq:a4}), ensuring we have accounted for all norm 4 vectors in $A_{24}$. Note that $a_2=2{25 \choose 2}=600$.\\

There are 36 vectors $u$ of type (i) satisfying $u\cdot w=-4$, namely $(-1^2,0^2,1^2,0^{19})$ and its permutations preserving both the group of coordinates 1 to 4, and that of coordinates 5 to 8. There are $34$ such vectors $u$ of type (ii) satisfying $u\cdot w=-4$: $17$ are obtained from $(-\nicefrac{4}{5}^4,\nicefrac{1}{5}^{20},-\nicefrac{4}{5})$ by allowing the last entry to be permuted to any of the last 17 coordinates, and $17$ are obtained from $(-\nicefrac{1}{5}^4,\nicefrac{4}{5}^{4},-\nicefrac{1}{5}^{16},\nicefrac{4}{5})$ in a similar fashion. Therefore $|S_4^w|=36+34=70$ and

\[
    \eta(L,w)=1+|S^w_2|+\frac{1}{2}|S_4^w|=1+16+\frac{1}{2}(70) = 52 \equiv 4 \not\equiv 0  \; (\text{mod 8}).
\]

\vspace{.2cm}

\noindent Thus $f_8(L)\geqslant w^2/2 = 4$. By Corollary \ref{cor:lowerbounds} (iii)-(iv) we conclude $g_4(A_{24})\geqslant 6$.\\

We may obtain the upper bound in two different ways, just as was done for $D_{24}$. First, as was mentioned in the introduction, $A_{24}$ is isomorphic to the canonical plumbed lattice associated to the Brieskorn sphere $-\Sigma(4,5,19)$, given as follows, with unmarked nodes of weight 2:
\begin{center}
\begin{tikzpicture}[scale=.5]

	\draw (0,0) -- (22,0);
	\draw (4,0) -- (4,1);
	
	\draw[fill=black] (0,0) circle(.1);
	\draw[fill=black] (1,0) circle(.1);
	\draw[fill=black] (2,0) circle(.1);
	\draw[fill=black] (3,0) circle(.1);
	\draw[fill=black] (4,0) circle(.1);
	\draw[fill=black] (5,0) circle(.1);
	\draw[fill=black] (6,0) circle(.1);
	\draw[fill=black] (7,0) circle(.1);
	\draw[fill=black] (8,0) circle(.1);
	\draw[fill=black] (9,0) circle(.1);
	\draw[fill=black] (10,0) circle(.1);
	\draw[fill=black] (11,0) circle(.1);
	\draw[fill=black] (12,0) circle(.1);
	\draw[fill=black] (13,0) circle(.1);
	\draw[fill=black] (14,0) circle(.1);
	\draw[fill=black] (15,0) circle(.1);
	\draw[fill=black] (16,0) circle(.1);
	\draw[fill=black] (17,0) circle(.1);
	\draw[fill=black] (18,0) circle(.1);
	\draw[fill=black] (19,0) circle(.1);
	\draw[fill=black] (20,0) circle(.1);
	\draw[fill=black] (21,0) circle(.1);
	\draw[fill=black] (22,0) circle(.1);
	\draw[fill=black] (4,1) circle(.1);
	
	\node at (3.5,1.25) {$4$};
	 
\end{tikzpicture}
\end{center}
Indeed, the node of weight 4 corresponds to $g$, and the nodes of weight 2 correspond to the roots $(0,1,-1,0^{22}),\ldots,(0^{23},1,-1)$. Then, as $-\Sigma(4,5,19)$ is diffeomorphic to $+1$ surgery on $T_{4,5}$ of slice genus 6, from inequality (\ref{eq:4gineq}) we obtain $g_4(A_{24})\leqslant 6$, implying equality.\\

For the second approach, begin with a smooth quintic in $\smash{\mathbb{C}\mathbb{P}^2}$. If only to relate this to our previous construction for $D_{24}$, we may start with 5 generic lines as in Figure \ref{fig:a24}, with no three lines intersecting in a common point, and then slightly perturb the defining equation to obtain a smooth quintic. Blowing up at 24 generic points yields a connected complex curve $\Sigma$ in $\smash{\mathbb{C}\mathbb{P}^2\# 24\overline{\mathbb{C}\mathbb{P}}^2}$ with

\[
    [\Sigma] \; = \; 5h - e_1-\cdots -e_{24} \; = \; \left(5 \, | \, 1^{24} \right).
\]

\vspace{.2cm}

\noindent By the adjunction formula (\ref{eq:adjunction}), the genus of $\Sigma$ is 6. The orthogonal complement of $[\Sigma]$ has roots $e_1-e_2,\ldots,e_{23}-e_{24}$ and $5h-2e_1-e_2-\cdots-e_{24}$, forming a copy of $-\mathsf{A}_{24}$. Being an even negative definite unimodular lattice, it must then be isomorphic to $-A_{24}$. From this construction we again conclude that $g_4(A_{24})\leqslant 6$. This completes the proof of Theorem \ref{thm:g4}, which, as discussed in the introduction, implies Theorem \ref{thm:2,11}.\\

\subsection{The lattice $A_{12}^2$}

Next we consider the Niemeier lattice $A_{12}^2$. We will show that $6\leqslant g_4(A_{12}^2)\leqslant 9$, thus providing another example of a Niemeier lattice not contained in $\mathscr{E}_{24}(-\Sigma(2,11,21))$.\\

The lattice $L=A_{12}^2$ is generated by its roots and the glue vector $g = \left((\nicefrac{1}{13}^{12}, -\nicefrac{12}{13}^1), (\nicefrac{5}{13}^8, -\nicefrac{8}{13}^5)\right)$. Note $g^2 = 4$. The norm 4 vectors are as follows, listed by type:

\begin{itemize}
    \item[(i)] $\phantom{\pm}((1^2, -1^2, 0^{9}),0)$
    \item[(ii)]$\phantom{\pm}(0,(1^2,-1^2,0^{9}))$
    \item[(iii)] $\phantom{\pm}((1,-1,0^{11}),(1,-1,0^{11}))$
    \item[(iv)] $\pm ((\nicefrac{1}{13}^{12}, -\nicefrac{12}{13}^1), (\nicefrac{5}{13}^8, -\nicefrac{8}{13}^5))$
    \item[(v)] $\pm((\nicefrac{2}{13}^{11}, -\nicefrac{11}{13}^2), (\nicefrac{10}{13}^3, -\nicefrac{3}{13}^{10}))$
    \item[(vi)] $\pm((\nicefrac{3}{13}^{10}, -\nicefrac{10}{13}^3), (\nicefrac{2}{13}^{11}, -\nicefrac{11}{13}^2))$
    \item[(vii)] $\pm((\nicefrac{5}{13}^8, -\nicefrac{8}{13}^5), (\nicefrac{12}{13}^1, -\nicefrac{1}{13}^{12}))$
\end{itemize}

\noindent Indeed, if we count the vectors listed, we compute

\begin{equation}
    \underbracket[.01cm]{\textstyle{13 \choose 4}{4 \choose 2}}_{\text{(i)}} +  \underbracket[.01cm]{\textstyle{13 \choose 4}{4 \choose 2}}_{\text{(ii)}} +  \underbracket[.01cm]{\textstyle2^2{13 \choose 2}^2}_{\text{(iii)}} +
    \underbracket[.01cm]{\textstyle 2(13){13 \choose 5}}_{\text{(iv)}} +
    \underbracket[.01cm]{\textstyle 2{13 \choose 2}{13 \choose 3}}_{\text{(v)}} +
    \underbracket[.01cm]{\textstyle 2{13 \choose 2}{13 \choose 3}}_{\text{(vi)}} + 
    \underbracket[.01cm]{\textstyle 2(13){13 \choose 5}}_{\text{(vii)}} \; = \; 189072.\label{eq:a122a4}
\end{equation}

\vspace{.2cm}

\noindent On the other hand, all roots in $A_{12}^2$ are of type $((1,-1,0^{11}),0)$ or $(0,(1,-1,0^{11}))$, and so $a_2=2{13 \choose 2}+ 2{13 \choose 2}=312$. From (\ref{eq:a4}) we compute $a_4=189072$, in agreement with (\ref{eq:a122a4}).\\

Consider $w=((\nicefrac{17}{13}, -\nicefrac{9}{13}^5, \nicefrac{4}{13}^7), (\nicefrac{7}{13}^6,-\nicefrac{6}{13}^7))$, of norm 8. We first check that $w\in A_{12}^2$. As seen in \cite[Ch.4 \S6]{conwaysloane}, $\mathsf{A}_n^\ast/\mathsf{A}_n$ is isomorphic to $\Z/(n+1)$, with a standard generator $[1]\in \mathbb{Z}/(n+1)$ corresponding to the class of a vector of type $(\nicefrac{1}{(n+1)}^n,-\nicefrac{n}{(n+1)}^1)$. Minimal norm vectors representing $[i]\in \Z/(n+1)$ are those of type $(\nicefrac{i}{(n+1)}^j,-\nicefrac{j}{(n+1)}^i)$ where $i+j=n+1$. Now, $A_{12}^2$ is by definition the preimage under $\oplus^2\mathsf{A}_{12}^\ast \to \oplus^2 \mathsf{A}_{12}^\ast/\mathsf{A}_{12}$ of the subgroup in the codomain generated by the class of $g$, which is $([1],[5])\in \oplus^2\Z/13$. This subgroup contains $([4],[7])$, so in particular $v:=((\nicefrac{4}{13}^2,-\nicefrac{9}{13}^{4},\nicefrac{4}{13}^7),(\nicefrac{7}{13}^6,-\nicefrac{6}{13}^7))\in A_{12}^2$. Finally, $w=((1,-1,0^{11}),0)+v$, so indeed $w\in A_{12}^2$.\\

It is straightforward to verify that $|w\cdot u| \leqslant u^2$ for all $u$ of norm 4 and norm 2, having listed all such vectors above. Thus $w$ is extremal by (\ref{eq:ex}). Further, the roots $u$ such that $u\cdot w=-2$ are obtained from $((-1,1,0^{11}),0)$ by permuting coordinates 2 to 6, and so $|S_2^w|=5$.\\

Next we count norm 4 vectors $u$ such that $w\cdot u=-4$ in order to compute $|S_4^w|$. The only contributions are from vectors of types (iv) and (vii); of the former are the 6 vectors obtained from $((-\nicefrac{12}{13},\nicefrac{1}{13}),(-\nicefrac{8}{13}^5,\nicefrac{5}{13}^8))$ by permuting the first 6 coordinates in the second $\mathsf{A}_{12}$-factor, and of the latter are the 6 vectors obtained from $((-\nicefrac{5}{13},\nicefrac{8}{13}^5,-\nicefrac{5}{13}^7),(-\nicefrac{12}{13},\nicefrac{1}{13}))$ using the same permutations in the second $\mathsf{A}_{12}$-factor. Thus $|S_4^w|=6+6=12$, and we compute

\[
    \eta(L,w) \; = \; 1 + |S^w_2| + \frac{1}{2}|S^w_4| \; = \; 1+5+\frac{1}{2}(12) \; = \; 12 \; \not\equiv \; 0 \; \text{ (mod 8)}.
\]

\vspace{.2cm}

\noindent Thus $f_8(A_{12}^2)\geqslant w^2/2 = 4$. By Corollary \ref{cor:lowerbounds} (iii)-(iv) we have $g_4(A_{12}^2)\geqslant 6$.\\

To obtain an upper bound, begin with $9$ lines in $\mathbb{C}\mathbb{P}^2$ as in the theorem of Pappus, with 9 triple points and some double points. Add a line passing through 2 double points and no triple points, and then add another line, now passing through one of the newly formed double points, and no other double or triple points. We end up with Figure \ref{fig:a122}. After blowing up triple points, and resolving double points, we obtain a connected complex curve $\Sigma$ in $\smash{\mathbb{C}\mathbb{P}^2\# 24\overline{\mathbb{C}\mathbb{P}}^2}$ with
\[
    [\Sigma] \; = \; 11h - 3\sum_{i=1}^{12}e_i -\sum_{i=13}^{24} e_i \; = \; \left(11 \, | \, 3^{12}, 1^{12} \right).
\]

\vspace{.0cm}

\noindent Using the adjunction formula (\ref{eq:adjunction}), we compute the genus of $\Sigma$ to be 9. The lattice orthogonal to $[\Sigma]$ is isomorphic to $-A_{12}^2$; the two copies of $-\mathsf{A}_{12}$ are spanned by $e_1-e_2,\ldots,e_{11}-e_{12}$, $(3 \, | \,1^{11}, 0^{13})$, and $e_{13}-e_{14},\ldots,e_{23}-e_{24}$, $(11 \, | \, 3^{12},2^1,1^{11})$,  respectively. Thus $g_4(A_{12}^2)\leqslant9$.

\subsection{The lattices $A_{15}D_9$ and $A_{17}E_7$}

\begin{figure}
\centering
\begin{minipage}{.5\textwidth}
  \centering
  \includegraphics[scale=.4]{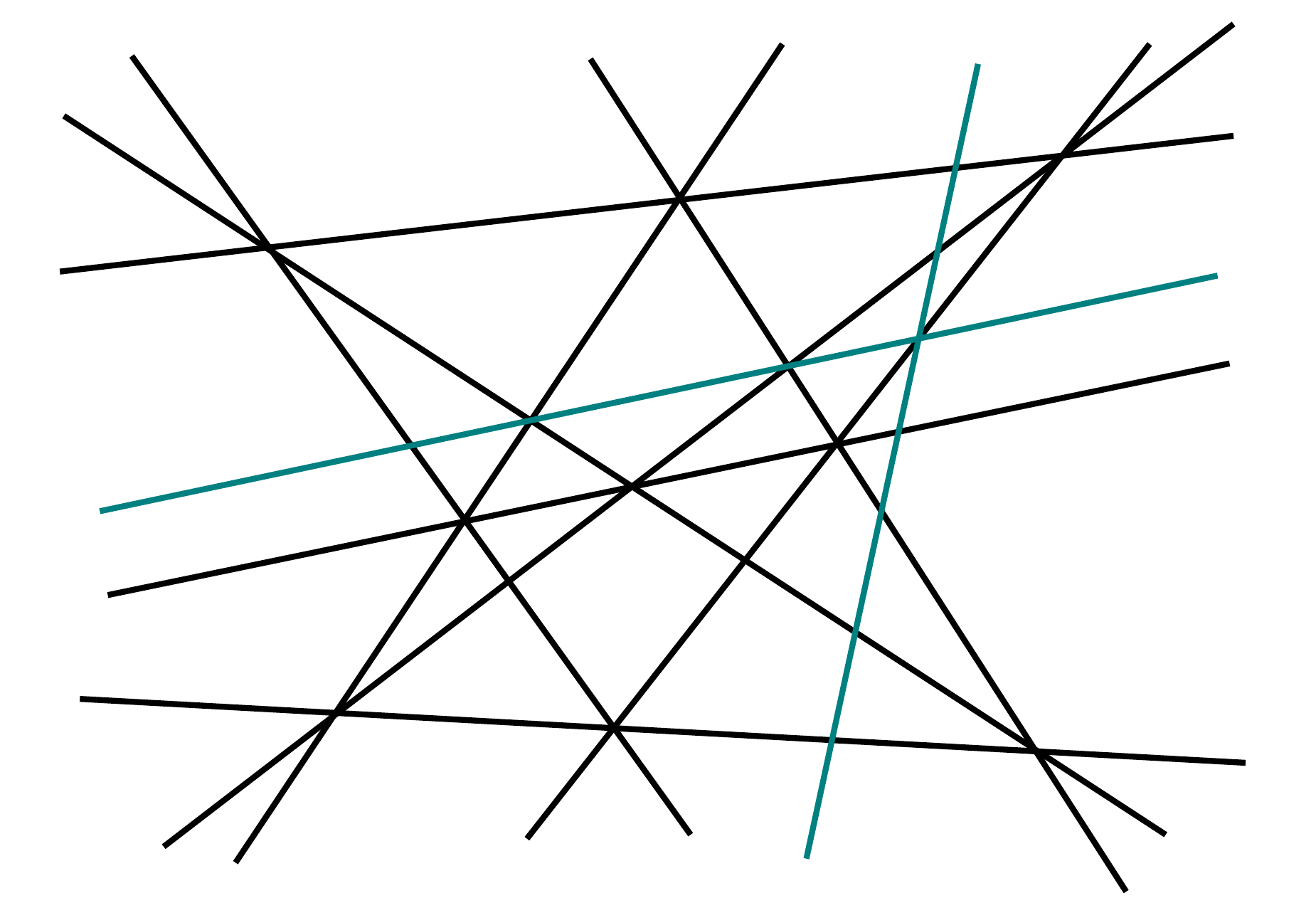}\vspace{.15cm}
  \captionof{figure}{$A_{12}^2$}
  \label{fig:a122}
\end{minipage}%
\begin{minipage}{.5\textwidth}
\centering
  \includegraphics[scale=.4]{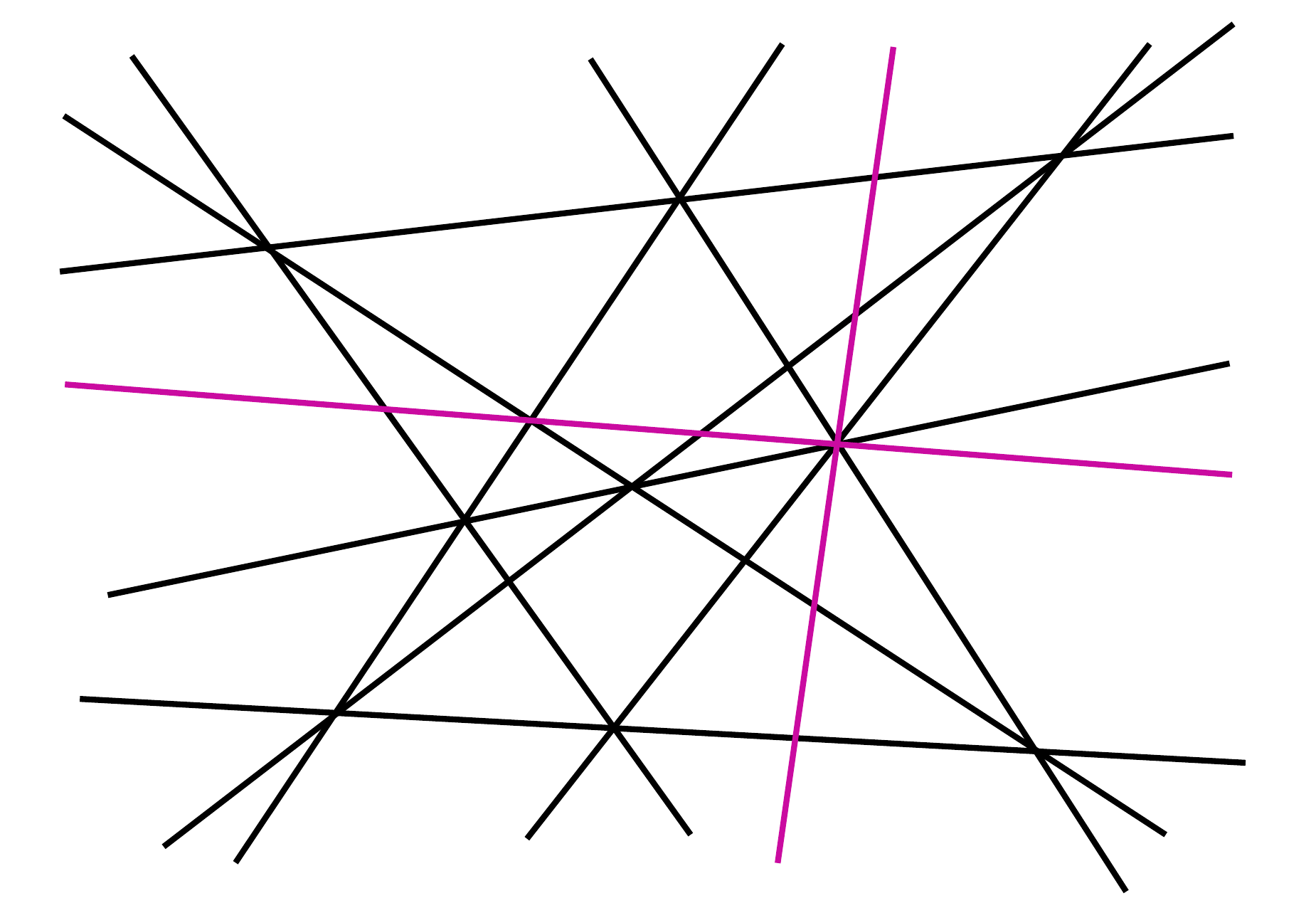}\vspace{.15cm}
  \captionof{figure}{$A_{15}D_9$}
  \label{fig:a15d9}
\end{minipage}
\end{figure}

We mention two upper bounds for $g_4(L)$ that fit into the family of constructions from above. First, consider the Niemeier lattice $A_{15}D_9$, generated by its roots and the glue vector $((\nicefrac{1}{8}^{14},-\nicefrac{7}{8}^2),(\nicefrac{1}{2}^9))$. Take a line arrangement with 11 lines, one point of multiplicity 5, 9 triple points, and some double points, as in Figure \ref{fig:a15d9}. Blow up all points not of multiplicity two and 14 generic points, and resolve double points, to obtain a complex curve in the class $(11 \, | \, 5, 3^9, 1^{14})$ of genus 8, by adjunction. The complement of this class is isomorphic to $-A_{15}D_9$, and thus $g_4(A_{15}D_9)\leqslant 8$.\\

\begin{figure}
\centering
\begin{minipage}{.5\textwidth}
  \centering
  \includegraphics[scale=.4]{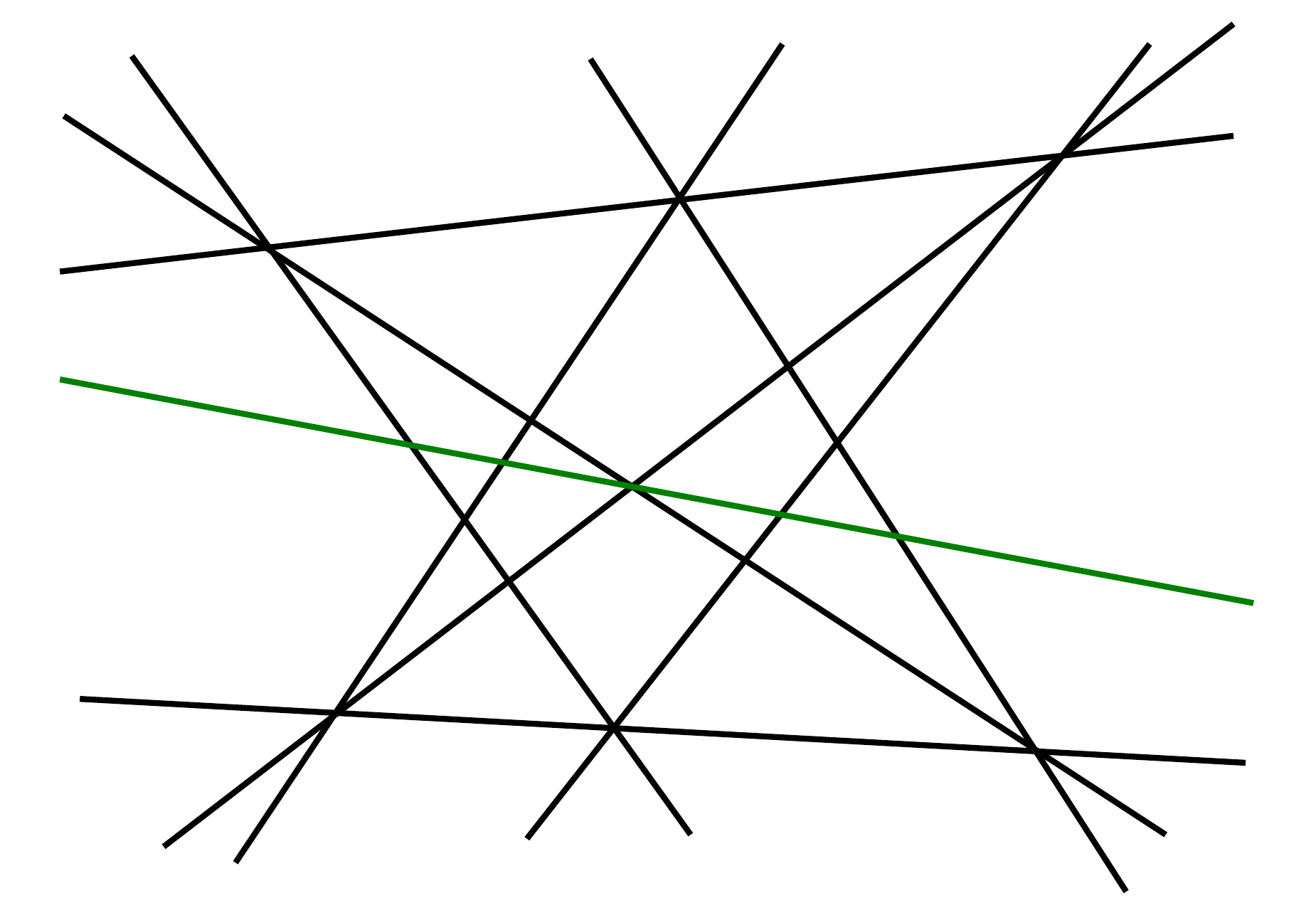}\vspace{.15cm}
  \captionof{figure}{$A_{17}E_7$}
  \label{fig:a17e7}
\end{minipage}%
\begin{minipage}{.5\textwidth}
\centering
  \includegraphics[scale=.4]{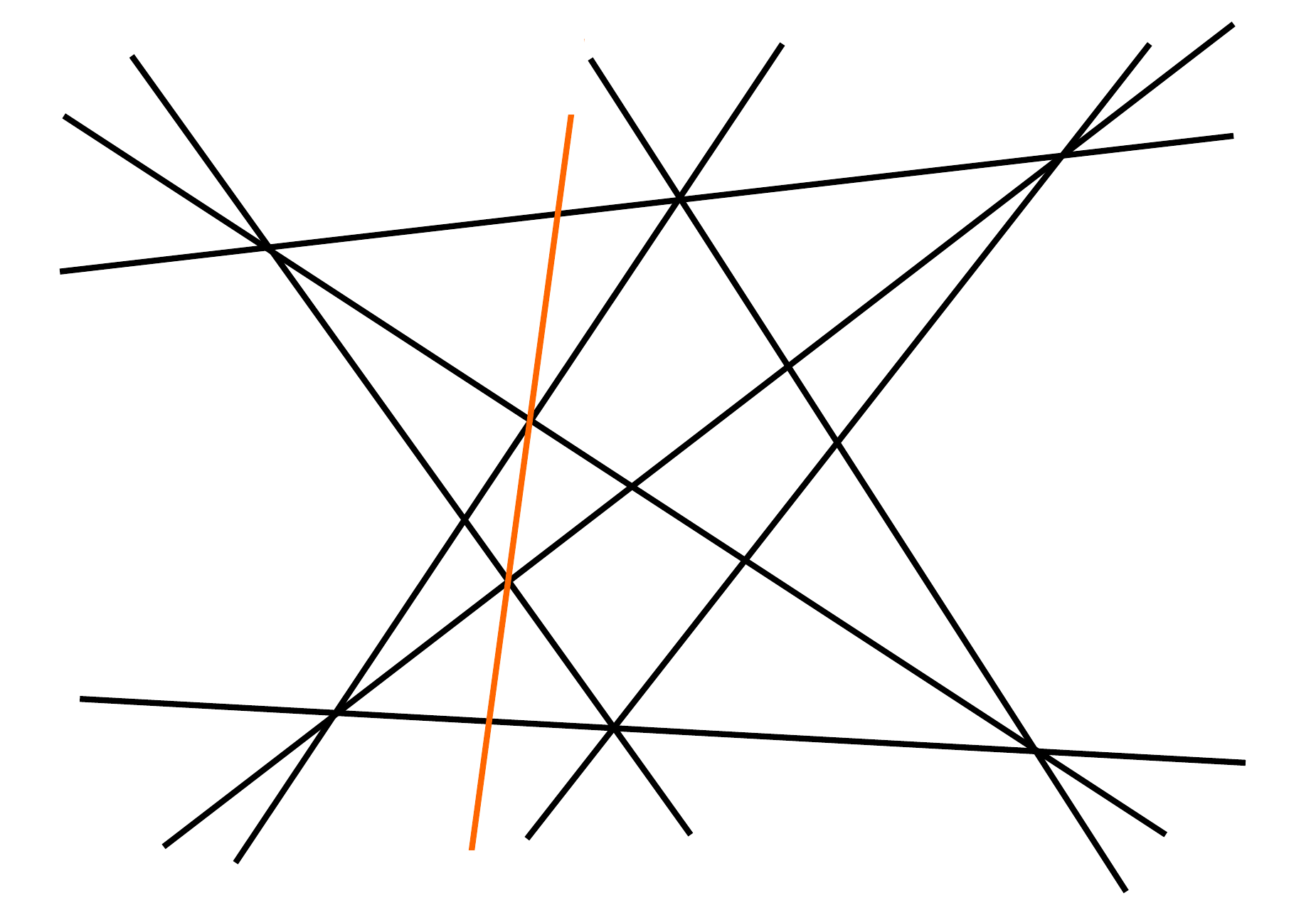}\vspace{.15cm}
  \captionof{figure}{$E_8^2$}
  \label{fig:e82}
\end{minipage}
\end{figure}

Similarly, we may consider the Niemeier lattice $A_{17}E_7$, generated by its roots and the glue vector $((\nicefrac{1}{6}^{15},-\nicefrac{5}{6}^3),(\nicefrac{1}{4}^6, -\nicefrac{3}{4}^2))$. Take an arrangement of 9 lines, interesecting 7 triple points and no other points of multiplicity greater than 2, as in Figure \ref{fig:a17e7}. Blow up the triple points and resolve the double points to obtain a curve in the class $(9\; | \; 3^7, 1^{17})$, of genus 7, by adjunction. The complement of this class is isomorphic to $-A_{17}E_7$, and thus $g_4(A_{17}E_7)\leqslant 7$. We do not know if either of these upper bounds are optimal. The Lorentzian forms for these lattices appeared in \cite[Ch.26]{conwaysloane}.

\subsection{The Leech lattice}

The Leech lattice $\Lambda_{24}$ is the unique positive definite unimodular lattice of rank 24 which is even and has no roots. We now show that the lower bound on $g_4(\Lambda_{24})$ obtained from our methods is no better than that from Heegaard Floer and Seiberg--Witten theory. The minimal norm of a nonzero vector in $\Lambda_{24}$ is $4$. Because $\Lambda_{24}$ is an even lattice, the next possible norms are $6$ and $8$.\\

\begin{theorem}[Ch.23 Thm.3\cite{conwaysloane}]\label{thm:conwaysloane} Every vector in $\Lambda_{24}$ is congruent mod $2\Lambda_{24}$ to a vector of norm $\leqslant 8$. The only congruences among such vectors are that each vector of norm $4$ or $6$ is congruent to its negative, and vectors of norm $8$ fall into congruence classes of size $48$, each class consisting of $24$ orthogonal vectors, up to signs.
\end{theorem}

\vspace{.25cm}

\noindent  This theorem implies that for $w\in\Lambda_{24}$ with $w^2=6$, we have $\text{Min}(w+2\Lambda_{24})=\{w,-w\}$, and thus $f_2(\Lambda_{24})\geqslant w^2-1 = 5$. In fact, the only larger extremal vectors $w$ have $\text{Min}(w+2\Lambda_{24})$ of size 48, and thus $f_2(\Lambda_{24})=5$. Similarly, we compute $f_4(\Lambda_{24})=f_8(\Lambda_{24})=6/2=3$. Then Corollary \ref{cor:lowerbounds} implies $g_4(\Lambda_{24})\geqslant 5$. However, as can be seen from the discussion in Section \ref{sec:oddsandends}, this is the same lower bound offered by Heegaard Floer and Seiberg--Witten theory. The author does not know of any upper bounds. We note that $-\Lambda_{24}$ is the orthogonal complement of $(145\; | \; 51, 47,45,\ldots,7,5,3)$, see \cite[Ch.26]{conwaysloane}.\\

\begin{remark} {\emph{Curiously, the lattice invariant $e(L)$ from \cite{froyshov-inequality} has $e(\Lambda_{24})=2$, as is easily computed from Theorem \ref{thm:conwaysloane}, and the instanton inequality of \cite[Thm.1]{froyshov-inequality} only provides the lower bound $g_4(\Lambda_{24})\geqslant 3$. This provides an example in which $e(L)\neq \delta(L)$, the latter being the lattice-theoretic quantity that appears in Heegaard Floer and Seiberg--Witten theory, reviewed in Section \ref{sec:oddsandends}.}}\\
\end{remark}

\subsection{Variations of Theorem \ref{thm:2,11}}\label{sec:otherknots}

We note that in proving Theorem \ref{thm:2,11}, the only properties of $T_{2,11}$ and $T_{4,5}$ used were that $T_{2,11}$ has slice genus 5 and $S_1^3(T_{2,11})$ bounds $D_{24}$, and that $S_1^3(T_{4,5})$ bounds $A_{24}$. To produce more examples, we may use constructions from \cite{golla-scaduto}.\\

For example, let $C$ be a rational quintic in $\mathbb{C}\mathbb{P}^2$ with cuspidal singularities whose links are knots $K_1,\ldots,K_n$. Blow up at 24 generic points of $C$ and remove a neighborhood of the proper transform to obtain a negative definite 4-manifold with boundary $-S_1^3(K_1\#\cdots\# K_n)$ and lattice $-A_{24}$; for more details see \cite{golla-scaduto}. Next, it is shown in \cite[Thm.2.3.10]{namba} that the possible singularity types $(K_1,\ldots,K_n)$ for a rational quintic are $(T_{4,5})$, $(T_{2,13})$, $(T_{3,5},T_{2,5})$, $(T_{3,4},T_{2,7})$, $(T_{2,9},T_{2,5})$, $(T_{3,4},T_{2,5},T_{2,3})$, $(T_{2,5},T_{2,5},T_{2,5})$, $(T_{2,7},T_{2,3},T_{2,3},T_{2,3})$. Note that in each of these cases, the knot $K_1\# \cdots \# K_n$ has slice genus 6, as follows, for example, from basic properties of Ozsv\'{ath}--Szab\'{o}'s $\tau$-invariant \cite{os-tau}, or Rasmussen's invariant \cite{rasmussen}.  \\

Similarly, if we begin with a rational septic in $\mathbb{C}\mathbb{P}^2$ with cusps of types $T_{5,6},K_1,\ldots,K_n$, we may blow up at the  $T_{5,6}$ cusp and 23 generic points to obtain a rational curve in the class $(7 \; | \; 5, 1^{23})$; in this way, we obtain as before a negative definite 4-manifold with boundary $-S_1^3(K_1\#\cdots\# K_n)$ and lattice $-D_{24}$. Examples are the rational septics $C_{7,1}$ and $C_{7,2}$ of \cite{flenner-zai}, which have $n=2$ and $(K_1,K_2)=(T_{2,3},T_{2,9})$ and $(K_1,K_2)=(T_{2,5},T_{2,7})$, respectively. Similarly to the remark at the end of the previous paragraph, the knot $K_1\# K_2$ in each case has slice genus 5. Write $mK$ for the $m$-fold connected sum of a knot $K$. We conclude:\\

\begin{corollary}
    Theorem \ref{thm:2,11} remains valid if $T_{4,5}$ is replaced by any of the knots
    \[
       T_{2,13}, \qquad T_{3,5}\# T_{2,5}, \quad T_{3,4}\# T_{2,7},\quad T_{2,9}\# T_{2,5}, \quad T_{3,4}\# T_{2,5}\# T_{2,3}, \quad 3T_{2,5}, \quad T_{2,7}\# 3 T_{2,3}
    \]
    and if $T_{2,11}$ is replaced by either $T_{2,3}\# T_{2,9}$ or $T_{2,5}\#T_{2,7}$.
\end{corollary}

\vspace{.2cm}

\vspace{.3cm}

\section{Odds and ends}\label{sec:oddsandends}

Although our focus in this article has been on Niemeier lattices, both the obstructions from instanton theory and the geometric constructions we have considered are applicable to all definite unimodular lattices. In this section we discuss some generalities regarding the geometric 4-genus $g_4(L)$ of minimal positive definite unimodular lattices $L$.\\

%First, we show that $g_4(L)$ is always finite. Let $L$ be a positive definite unimodular lattice. From the classification of indefinite lattices, $L\oplus \langle -1\rangle$ is isomorphic to $\langle 1\rangle^n\oplus\langle -1\rangle$ where $n=\text{rk}L$. Under this isomorphism, $L$ is the orthogonal complement of a vector $\langle 1\rangle^n\oplus\langle -1\rangle$ with $v^2=-1$. Now we may consider $v$ to be an element in the lattice of $\smash{-\mathbb{C}\mathbb{P}^2\# n\overline{\mathbb{C}\mathbb{P}}^2}$. Choose a surface representative $\Sigma$ of $v$; then $g_4(L)\leqslant \text{genus}(\Sigma)$.\\

For completeness, we record the lower bound imposed on $g_4(L)$ that comes from Seiberg-Witten and Heegaard Floer theory. Recall that a {\emph{characteristic vector}} $\xi\in L$ is a vector such that $\xi\cdot x \equiv x^2$ (mod 2) for all $x\in L$. Denote by $\text{Char}(L)$ the set of characteristic vectors for $L$. Define

\[
    \delta(L) \;:= \; \max_{\xi\in \text{Char}(L)}\frac{1}{8}(\text{rk} L - \xi^2).
\]

\vspace{.2cm}

\noindent It turns out that $\delta(L)$ is a non-negative integer. A result of Elkies \cite{elkies-1} says that $\delta(L)=0$ if and only if $L$ is diagonal. If $L$ is even, we note that $\delta(L)=\text{rk}L/8$. We then have the inequality

\begin{equation}
    g_4(L) \;\geqslant\; 2\delta(L)-1 \label{eq:dinvineq}
\end{equation}

\vspace{.2cm}

\noindent as follows from \cite[Thm.9.15]{os} and \cite[Prop.3.4]{rasmussen}; see also \cite[Thm.1.1,6.1]{bg}, which applies after blowing up our surface $\Sigma$ once to obtain a surface of self-intersection 0. In particular, for an even lattice $L$, we have $g_4(L)\geqslant \text{rk} L /4 -1$.\\

Perhaps the most general question in the current context is the following: what range of values $(\text{rk} L, g_4(L))$ occur as $L$ runs through all minimal positive definite unimodular lattices? What is the range when $L$ is restricted to even lattices? These are \textbf{{\color{blue}geography}} problems, to which Theorem \ref{thm:g4} provides some very small progress when the rank is fixed at 24, distinct from the information provided by Seiberg-Witten and Heegaard Floer theory above.\\

There are also \textbf{{\color{blue}botany}} problems. For example, we may ask: for fixed $g\in\Z_{\geqslant 0}$, how many minimal positive definite unimodular lattices have $g_4(L)=g$? For low values of $g$, this problem is solved. Indeed, for a minimal positive definite unimodular lattice $L$, we have
\begin{align}
    g_4(L) \; = \; 0 \qquad & \Longleftrightarrow \qquad L= 0 \label{eq:g40}\\
    g_4(L) \; = \; 1 \qquad  & \Longleftrightarrow \qquad L=E_8 \label{eq:g41}\\
    g_4(L) \; = \; 2 \qquad & \Longleftrightarrow \qquad L = \Gamma_{12} \label{eq:g42}
\end{align}

\vspace{.1cm}

\noindent  Here $\Gamma_{4k}$ is the lattice generated by $\mathsf{D}_{4k}$ and the glue vector $(\nicefrac{1}{2}^{4k})$. Note that $\Gamma_{4k}$ is even precisely when $k$ is even. The equivalence (\ref{eq:g40}) follows from Donaldson's theorem, cf. \cite[Prop.2.2.11]{gs}; alternatively, we may use (\ref{eq:dinvineq}) and the above cited result of Elkies. That of (\ref{eq:g41}) is related to the main result of Fr\o yshov's thesis \cite{froyshov-thesis}, and follows from (\ref{eq:mainineq2}), Corollary \ref{cor:lowerbounds} and Lemmas 3.2 and 3.3 of \cite{scaduto-forms}. The equivalence (\ref{eq:g42}) is related to the main result of \cite{scaduto-forms}, and follows from (\ref{eq:mainineq}) with $n=4$, Corollary \ref{cor:lowerbounds} and Lemma 4.2 of \cite{scaduto-forms}. \\

We next ask for positive definite lattices $L$ with $g_4(L)=3$. We do not solve this problem here, but provide several examples, some of which have $g_4(L)=3$, and some of which we only know to have $g_4(L)\in\{3,4\}$. We also include some examples with $g_4(L)=4$, and one with $g_4(L)\in\{4,5\}$.\\

\subsubsection{The lattices $\Gamma_{16}$ and $E_8^2$}

First consider the case of even lattices. There are two such lattices of rank 16: $\Gamma_{16}$ and $E_8^2$. We have $g_4(\Gamma_{16})=3$, as follows, for example, from \cite[Prop.4.1]{scaduto-forms}. For the other rank 16 even unimodular lattice we have $g_4(E_8^2)\in \{3,4\}$. Indeed, for the lower bound, we may either use (\ref{eq:dinvineq}), or compute $f_4(E_8^2)= 2$ and use Corollary \ref{cor:lowerbounds}. For the upper bound, we may take 9 lines in $\mathbb{C}\mathbb{P}^2$ with 8 triple points and some double points, as in Figure \ref{fig:e82}. After blowing up at the triple points, 8 generic points, and resolving double points, we obtain a complex curve $\Sigma$ in $\smash{\mathbb{C}\mathbb{P}^2\# 16\overline{\mathbb{C}\mathbb{P}}^2}$ with

\[
    [\Sigma] \; = \; 9h - 3\sum_{i=1}^8 e_i -\sum_{i=9}^{16}e_i \; = \; ( 9\, | \, 3^8, 1^8).
\]

\vspace{.2cm}

\noindent The complement lattice of $[\Sigma]$ is isomorphic to $-E_8^2$. Indeed, one copy of $-E_8$ is spanned by $e_1-e_2,\ldots,e_7-e_8,h-e_1-e_2-e_3$, and the other by $e_9-e_{10},\ldots,e_{15}-e_{16}, (3 \, | \, 1^{11}, 0^{5})$. Furthermore, the adjunction formula implies that $\Sigma$ has genus 4. Thus $g_4(E_8^2)\leqslant 4$.\\

\begin{figure}
\centering
\begin{minipage}{.5\textwidth}
  \centering
  \includegraphics[scale=.4]{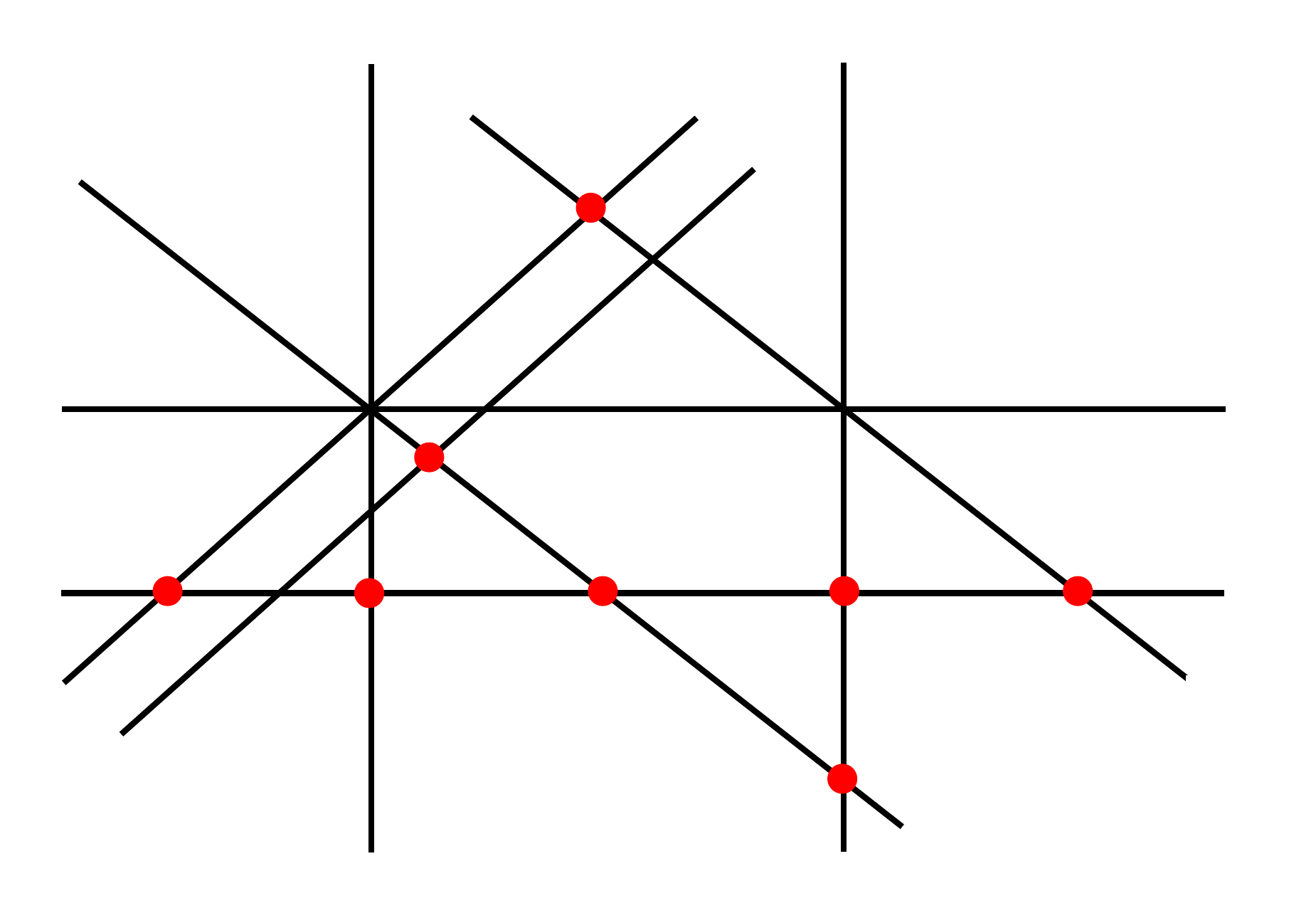}\vspace{.15cm}
  \captionof{figure}{$D_8^2$}
  \label{fig:d82}
\end{minipage}%
\begin{minipage}{.5\textwidth}
\centering
  \includegraphics[scale=.4]{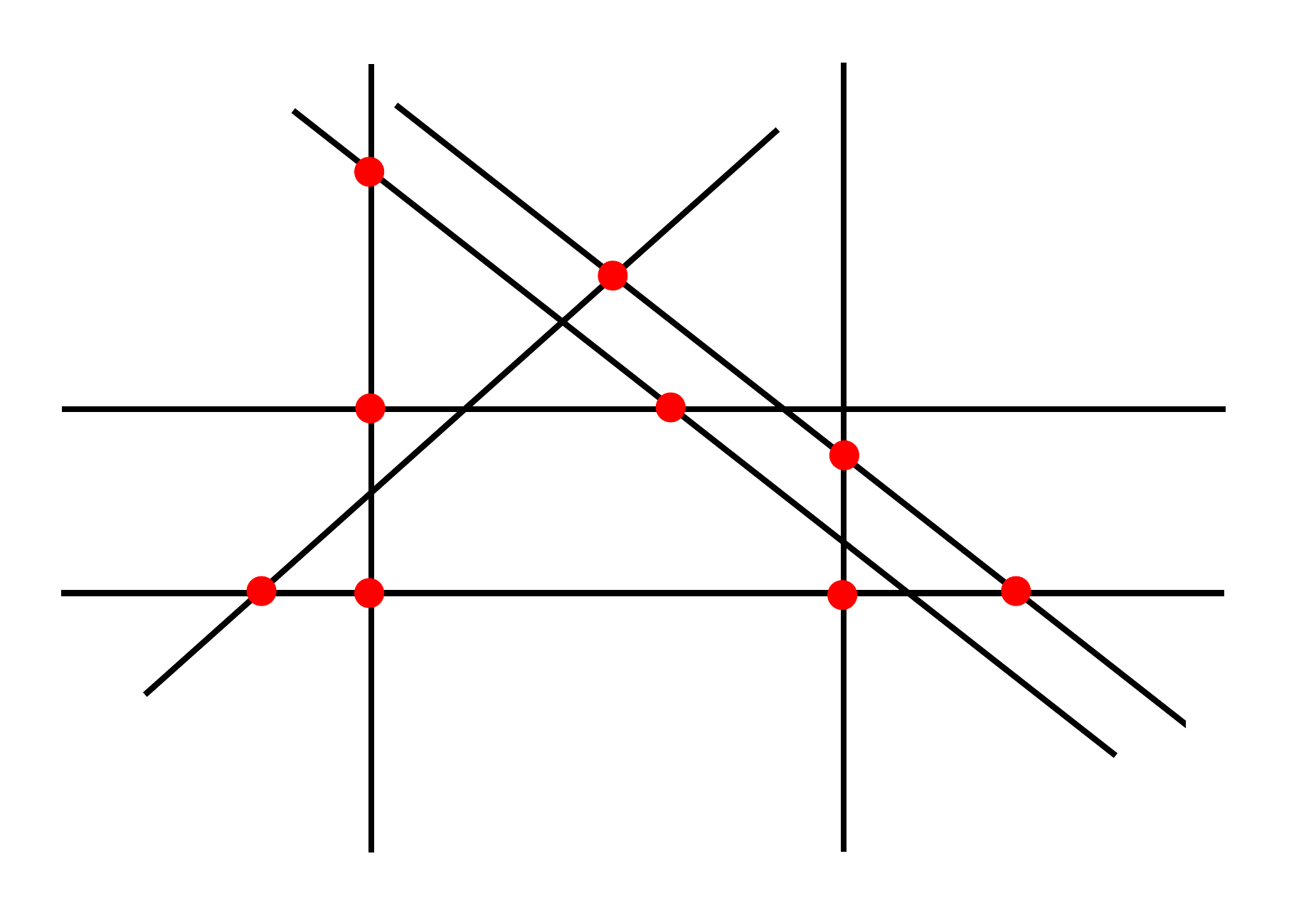}\vspace{.15cm}
  \captionof{figure}{$A_9^2$}
  \label{fig:a92}
\end{minipage}
\end{figure}

\subsubsection{The lattice $A_{15}$} This is the unique positive definite unimodular lattice with root lattice $\mathsf{A}_{15}$, an example considered in \cite[\S 5]{scaduto-forms}. It is generated by its roots and the glue vector $(\nicefrac{1}{4}^{12},-\nicefrac{3}{4}^{4})$, and is isomorphic to the canonical plumbing for $-\Sigma(3,4,11)$, which is $+1$ surgery on $T_{3,4}$. As $g_4(T_{3,4})=3$, from (\ref{eq:4gineq}) we have $g(A_{15})\leqslant 3$. Alternatively, we may start with a quartic in $\mathbb{C}\mathbb{P}^2$ and blow up at 15 generic points to obtain a complex curve of genus 3 representing the class $(4 \, | \, 1^{15})$, whose complementary lattice is isomorphic to $-A_{15}$. The lower bound follows from $f_4(A_{15})=2$ and Corollary \ref{cor:lowerbounds}.\\

\subsubsection{The lattice $E_7^2$} The next example is provided by $E_{7}^2$, the positive definite unimodular lattice with full-rank root lattice $\mathsf{E}_7\oplus \mathsf{E}_7$. It is generated by its roots and the glue vector $((\nicefrac{1}{4}^6,-\nicefrac{3}{4}^2),(\nicefrac{1}{4}^6,-\nicefrac{3}{4}^2))$. To obtain an upper bound for $g_4$, begin with a quartic and conic in $\mathbb{C}\mathbb{P}^2$ in general position, intersecting in 8 points. Blow up 7 of these double points, and resolve the last, to obtain a connected complex curve representing the class $(6 \, | \, 2^7, 1^7)$, of genus 3. The complement lattice is isomorphic to $-E_7^2$. Indeed, one copy of $-\mathsf{E}_7$ is spanned by the roots $e_1-e_2,\ldots, e_{6}-e_7,(1 \, | \, 1^{3},0^{11})$, and the other by $e_8-e_9,\ldots,e_{13}-e_{14},(3 \, | \, 1^{11},0^3)$. The lower bound $g_4(E_7^2)\geqslant 3$ is provided by $f_4(E_7^2)=2$, see \cite[Prop.9.1]{scaduto-forms}, and Corollary \ref{cor:lowerbounds}.\\

The upper bound $g_4(E_7^2)\leqslant 3$ may also be obtained from (\ref{eq:4gineq}) and the fact that $E_7^2\in \mathscr{L}(S_1^3(K))$ for several knots $K$ of slice genus 3, for example, $K=T_{3,4}$; see \cite[\S 4.3]{golla-scaduto}.\\

Finally, we note that $E_7^2$ is isomorphic to the following plumbed lattice:
\begin{center}
\begin{tikzpicture}[scale=0.65]
	\draw (0,0) -- (1,1);
	\draw (0,0) -- (1,0);
	\draw (0,0) -- (1,-1);
	\draw (1,1) -- (3,1);
	\draw (1,0) -- (4,0);
	\draw (1,-1) -- (6,-1);
	
	\draw[fill=black] (0,0) circle(.1);
	\draw[fill=black] (1,0) circle(.1);
	\draw[fill=black] (1,1) circle(.1);
	\draw[fill=black] (3,1) circle(.1);
	\draw[fill=black] (1,-1) circle(.1);
	\draw[fill=black] (2,0) circle(.1);
	\draw[fill=black] (2,1) circle(.1);
	\draw[fill=black] (3,0) circle(.1);
	\draw[fill=black] (4,0) circle(.1);
	\draw[fill=black] (2,-1) circle(.1);
	\draw[fill=black] (3,-1) circle(.1);
	\draw[fill=black] (4,-1) circle(.1);
	\draw[fill=black] (5,-1) circle(.1);
	\draw[fill=black] (6,-1) circle(.1);
	
	\node[font=\small] at (-.45,0) {};
	\node[font=\small] at (1,1.4) {};
	\node[font=\small] at (2,1.4) {$3$};
	\node[font=\small] at (3,0.4) {$3$};
	\node[font=\small] at (2,-.6) {$3$};
	
	\node at (-.5,1) {$E_7^2$};
\end{tikzpicture}
\end{center}

\vspace{.2cm}

\noindent The unmarked nodes are of weight 2. Indeed, a basis for this plumbing is given by
\begin{gather*}
 ((\nicefrac{1}{4}^2,-\nicefrac{3}{4},\nicefrac{1}{4}^4,-\nicefrac{3}{4}),(-\nicefrac{3}{4},\nicefrac{1}{4}^5,-\nicefrac{3}{4},\nicefrac{1}{4})),\\
 ((\nicefrac{3}{4},-\nicefrac{1}{4}^5,\nicefrac{3}{4},-\nicefrac{1}{4}),(-\nicefrac{1}{4}^6,\nicefrac{3}{4}^2)), \qquad
((\nicefrac{3}{4}^2,-\nicefrac{1}{4}^6),(\nicefrac{3}{4}^2,-\nicefrac{1}{4}^6)),
\end{gather*}

\vspace{.1cm}

\noindent which account for the nodes of weight 3, along with the roots $-e_i+e_{i+1}$ for $i\in\{1,10,11,12,13\}$, $e_i-e_{i+1}$ for $i\in\{3,4,5,6,15\}$, and $((-\nicefrac{1}{2}^4,\nicefrac{1}{2}^4),0)$. Here we view $E_7^2$ as embedded in $\R^{16}$, and we write $e_1,\ldots,e_{16}$ for the standard basis of $\R^{16}$. On the other hand, the plumbing given is the canonical positive definite plumbing for $-\Sigma(8,11,17)$, and so we have $-E_7^2\in \mathscr{L}(\Sigma(8,11,17))$.\\

\subsubsection{The lattice $D_8^2$} This lattice is generated by its roots and the glue vectors $((\nicefrac{1}{2}^8),(1,0^7))$ and $((1,0^7),(\nicefrac{1}{2}^8))$. As $f_4(D_8^2)\geqslant 2$, see \cite[Lemma 4.2]{scaduto-forms}, we have $g_4(D_8^2)\geqslant 3$ from Corollary \ref{cor:lowerbounds}. To obtain an upper bound, take 8 lines in $\mathbb{C}\mathbb{P}^2$ with one quadruple point, one triple point, and all other intersections of multiplicity 2, as in Figure \ref{fig:d82}. Blow up the quadruple and triple points, the 8 double points indicated in Figure \ref{fig:d82}, and 6 generic points, and resolve the remaining double points to obtain a connected complex curve $\Sigma$ representing $(8 \, |  \, 4,3,2^8,1^6)$. (Any selection of double points such that the resulting curve is connected is fine.) By adjunction, the genus of $\Sigma$ is 4. The upper bound $g_4(D_8^2)\leqslant 4$ then follows from:\\

\begin{lemma}
    The complement of $v=(8\; | \; 4,3,2^8, 1^6)\in H_2(\mathbb{C}\mathbb{P}^2\# 16 \overline{\mathbb{C}\mathbb{P}}^2;\Z)$ is isomorphic to $-D_8^2$.
\end{lemma}

\vspace{.1cm}

\begin{proof}
First note that $v^2=1$; thus the complement is unimodular and negative definite. Next, note that the roots $e_3-e_4,\ldots,e_9-e_{10},h-e_1-e_3-e_4$ span one copy of $-\mathsf{D}_8$, and $e_{11}-e_{12},\ldots, e_{15}-e_{16}, h-e_1-e_2-e_{11}, (3\, | \, 1^{11}, 0^5),(4\, | \, 2,1^{14},0)$ span another, orthogonal copy. Now there are four negative definite unimodular lattices containing two orthogonal copies of $-\mathsf{D}_8$, namely $-D_8^2$, $-E_8^2$, $-E_8\oplus \langle -1 \rangle^8$ and $\langle -1 \rangle^{16}$. The latter two possibilities only occur if one half the sum of two short legs in any root basis for one copy of $-\mathsf{D}_8$ lies in $L$. For one copy we may take these short legs to be $e_3-e_4$ and $h-e_1-e_3-e_4$; for the other, $h-e_1-e_3-e_{11}$ and $3h-e_1-\cdots-e_{11}$. Half the sum of each pair is not contained in $L$. Thus our lattice is either $-D_8^2$ or $-E_8^2$. However, $h-e_3-e_4-e_5-e_6\in L$ is of norm 3, so $L$ is an odd lattice, and must be isomorphic to $-D_8^2$, as claimed.
\end{proof}

\vspace{.2cm}

We also note that the upper bound $g_4(D_8^2)\leqslant 4$ may also be obtained using (\ref{eq:4gineq}) and that fact that $D_8^2\in \mathscr{L}(S_1^3(K))$ for several knots $K$ of slice genus 4, such as $K=T_{3,5}$; see \cite[\S 4.3]{golla-scaduto}.\\

As in the case of $E_7^2$, the lattice $D_8^2$ may be represented by the following plumbing:
\begin{center}
\begin{tikzpicture}[scale=0.65]
	\draw (0,0) -- (1,1);
	\draw (0,0) -- (1,0);
	\draw (0,0) -- (1,-1);
	\draw (1,1) -- (3,1);
	\draw (1,0) -- (4,0);
	\draw (1,-1) -- (8,-1);
	
	\draw[fill=black] (0,0) circle(.1);
	\draw[fill=black] (1,0) circle(.1);
	\draw[fill=black] (1,1) circle(.1);
	\draw[fill=black] (3,1) circle(.1);
	\draw[fill=black] (1,-1) circle(.1);
	\draw[fill=black] (2,0) circle(.1);
	\draw[fill=black] (2,1) circle(.1);
	\draw[fill=black] (3,0) circle(.1);
	\draw[fill=black] (4,0) circle(.1);
	\draw[fill=black] (2,-1) circle(.1);
	\draw[fill=black] (3,-1) circle(.1);
	\draw[fill=black] (4,-1) circle(.1);
	\draw[fill=black] (5,-1) circle(.1);
	\draw[fill=black] (6,-1) circle(.1);
	\draw[fill=black] (7,-1) circle(.1);
	\draw[fill=black] (8,-1) circle(.1);
	
	\node[font=\small] at (-.45,0) {};
	\node[font=\small] at (1,1.4) {};
	\node[font=\small] at (2,1.4) {};
	\node[font=\small] at (4,0.4) {$3$};
	\node[font=\small] at (1,-.6) {$3$};
	
	\node at (-.5,1) {$D_8^2$};
\end{tikzpicture}
\end{center}

\vspace{.1cm}

\noindent Indeed, a basis for this plumbing is provided by the square 3 vectors $((-1,0^7),(-\nicefrac{1}{2}^4,\nicefrac{1}{2}^4))$ and $((\nicefrac{1}{2}^8),(-1,0^7))$, and the roots $e_i-e_{i+1}$ for $1\leqslant i \leqslant 7$ and $9\leqslant i \leqslant 15$. Here we view $D_8^2$ as embedded in $\R^{16}$. This is the canonical positive definite plumbing for $-\Sigma(4,9,17)$. Thus $-D_8^2\in \mathscr{L}(\Sigma(4,9,17))$.\\

\subsubsection{The lattice $A_{17}A_1$}

The Brieskorn sphere $-\Sigma(3,5,14)$ has canonical positive definite lattice the following plumbing:
\begin{center}
\begin{tikzpicture}[scale=0.65]
	\draw (0,0) -- (1,1);
	\draw (0,0) -- (1,0);
	\draw (0,0) -- (1,-1);
	\draw (1,1) -- (2,1);
	\draw (1,0) -- (2,0);
	\draw (1,-1) -- (13,-1);
	
	\draw[fill=black] (0,0) circle(.1);
	\draw[fill=black] (1,0) circle(.1);
	\draw[fill=black] (1,1) circle(.1);
	\draw[fill=black] (1,-1) circle(.1);
	\draw[fill=black] (2,0) circle(.1);
	\draw[fill=black] (2,1) circle(.1);
	\draw[fill=black] (2,-1) circle(.1);
	\draw[fill=black] (3,-1) circle(.1);
	\draw[fill=black] (4,-1) circle(.1);
	\draw[fill=black] (5,-1) circle(.1);
	\draw[fill=black] (6,-1) circle(.1);
	\draw[fill=black] (7,-1) circle(.1);
	\draw[fill=black] (8,-1) circle(.1);
	\draw[fill=black] (9,-1) circle(.1);
	\draw[fill=black] (10,-1) circle(.1);
	\draw[fill=black] (11,-1) circle(.1);
	\draw[fill=black] (12,-1) circle(.1);
	\draw[fill=black] (13,-1) circle(.1);
	
	\node[font=\small] at (-.45,0) {};
	\node[font=\small] at (1,1.4) {};
	\node[font=\small] at (2,1.4) {};
	\node[font=\small] at (1,0.4) {$3$};
	
	\node at (-.5,1) {$A_{17}A_1$};
\end{tikzpicture}
\end{center}

\vspace{.2cm}

\noindent The roots in this plumbing together span $\mathsf{A}_{16}\oplus \mathsf{A}_1$. As can be seen from \cite[p.416]{conwaysloane}, the only positive definite unimodular lattice containing this root lattice is $L=A_{17}A_1$. This lattice is generated by its roots and the glue vector $((\nicefrac{1}{6}^{15},-\nicefrac{5}{6}^3),(\nicefrac{1}{2},-\nicefrac{1}{2}))$, of norm 3. We claim $g_4(L)= 4$.\\

To establish the lower bound, consider $w=((\nicefrac{7}{6}, \nicefrac{1}{6}^{13},-\nicefrac{5}{6}^4),(\nicefrac{1}{2},-\nicefrac{1}{2}))$, of norm 5. Then $w$ is extremal and the roots $u$ such that $u\cdot w=-2$ are $((-1,0^{13},1,0^3),0)$, $((-1,0^{14},1,0^2),0)$, $((-1,0^{15},1,0),0)$, $((-1,0^{16},1),0)$. Thus $|S_2^w|=4$, and $|\text{Min}(L+2w)/\pm |\equiv 1+|S_2^w| \equiv 1 + 4 \equiv 1$ (mod 2). Then $f_2(L)\geqslant w^2-1 = 4$, and by Corollary \ref{cor:lowerbounds} (i), we have $g_4(A_{17}A_1)\geqslant 4$.\\

On the other hand, $-\Sigma(3,5,14)$ is $+1$ surgery on the torus knot $T_{3,5}$ of slice genus 4, and so by implication \eqref{eq:4gineq}, we have $g_4(A_{17}A_1)\leqslant 4$, implying equality.\\

\subsubsection{The lattice $A_9^2$}

The lattice $L=A_9^2$ is the rank 18 positive definite unimodular lattice generated by $\mathsf{A}_9\oplus \mathsf{A}_9$ and the glue vector $g=((-\nicefrac{7}{10}^3, \nicefrac{3}{10}^7),(-\nicefrac{9}{10},\nicefrac{1}{10}^9))$. Thus $A_9^2$ is the preimage under $\oplus^2\mathsf{A}^\ast_9\to \oplus^2 \mathsf{A}^\ast_9/\mathsf{A}_9$ of the subgroup in the codomain generated by $x=([3],[1])$. A minimal norm representative for $([5],[5])=5x$ is $w=((\nicefrac{1}{2}^5,-\nicefrac{1}{2}^5),(\nicefrac{1}{2}^5,-\nicefrac{1}{2}^5))$, of norm 5. We see that $S_2^w=\emptyset$, and so $\text{Min}(w+2L)=\{w,-w\}$, implying $f_2(L)\geqslant w^2-1=4$. By Corollary \ref{cor:lowerbounds} (i) we obtain $g_4(L)\geqslant 4$.\\ 

To obtain an upper bound, choose 7 generic lines in $\mathbb{C}\mathbb{P}^2$ having only double points, as in Figure \ref{fig:a92}. Blow up at the 10 indicated double points, resolve the remaining ones, and blow up at 8 generic points to obtain a connected complex curve $\Sigma$ in $\smash{\mathbb{C}\mathbb{P}^2\# 18\overline{\mathbb{C}\mathbb{P}}^2}$ representing the class $(7\;|\; 2^{10},1^8)$. The complement of this class has one copy of $\mathsf{A}_9$ spanned by $e_1-e_2,\ldots,e_9-e_{10}$ and another, orthogonal copy spanned by $e_{11}-e_{12},\ldots, e_{17}-e_{18},(3\; | \; 1^{11},0^7)$. Thus $[\Sigma]^\perp$ is isomorphic to $A_9^2$, the unique positive definite unimodular lattice containing this root lattice, see e.g. \cite[p.416]{conwaysloane}. By adjunction, the genus of $\Sigma$ is 5, and so $g_4(A_9^2)\leqslant 5$.\\

We mention that $A_9^2$ is isomorphic to the following plumbed lattice:
\begin{center}
\begin{tikzpicture}[scale=0.65]
	\draw (0,0) -- (1,1);
	\draw (0,0) -- (1,0);
	\draw (0,0) -- (1,-1);
	\draw (1,1) -- (2,1);
	\draw (1,0) -- (6,0);
	\draw (1,-1) -- (9,-1);
	
	\draw[fill=black] (0,0) circle(.1);
	\draw[fill=black] (1,0) circle(.1);
	\draw[fill=black] (1,1) circle(.1);
	\draw[fill=black] (1,-1) circle(.1);
	\draw[fill=black] (2,0) circle(.1);
	\draw[fill=black] (2,1) circle(.1);
	\draw[fill=black] (3,0) circle(.1);
	\draw[fill=black] (4,0) circle(.1);
	\draw[fill=black] (5,0) circle(.1);
	\draw[fill=black] (6,0) circle(.1);
	\draw[fill=black] (2,-1) circle(.1);
	\draw[fill=black] (3,-1) circle(.1);
	\draw[fill=black] (4,-1) circle(.1);
	\draw[fill=black] (5,-1) circle(.1);
	\draw[fill=black] (6,-1) circle(.1);
	\draw[fill=black] (7,-1) circle(.1);
	\draw[fill=black] (8,-1) circle(.1);
	\draw[fill=black] (9,-1) circle(.1);
	
	\node[font=\small] at (-.45,0) {};
	\node[font=\small] at (1,1.4) {};
	\node[font=\small] at (2,1.4) {};
	\node[font=\small] at (1,-0.6) {$3$};
	
	\node at (-.5,1) {$A_9^2$};
\end{tikzpicture}
\end{center}
A corresponding basis is given by $e_i-e_{i+1}$ for $1\leqslant i \leqslant 9$, $10\leqslant i\leqslant 16$, and the glue vector $g$, which corresponds to the node of weight 3. This is the canonical positive definite plumbing for $-\Sigma(3,7,19)$, so we conclude that $-A_9^2\in \mathscr{L}(\Sigma(3,7,19))$.\\

The (in)equalities obtained above are summarized in the following proposition. We include the rank 20 lattice $\Gamma_{20}$, which has $g_4=4$, as computed in \cite[Prop.4.1]{scaduto-forms}.\\

\begin{prop} The lattices $\Gamma_{16}$, $E_7^2$ and $A_{15}$ have $g_4=3$. The lattices $E_8^2$ and $D_8^2$ have $g_4\in\{3,4\}$. The lattices $A_{17}A_1$ and $\Gamma_{20}$ have $g_4=4$. The lattice $A_9^2$ has $g_4\in \{4,5\}$. \\
\end{prop}

%\newpage

%\section{Scraps}\label{sec:scraps}

%\input{scraps.tex}

\bibliographystyle{alpha}

\end{document}